\documentclass[12pt]{article}

\usepackage{tabularx}
\usepackage{array}
\usepackage{booktabs}

\usepackage{amsmath, amssymb}
\usepackage[dvips]{graphicx}
\usepackage{amsthm, amscd}

\usepackage{bm}

\newcommand{\calA}{\mathcal{A}}
\newcommand{\calB}{\mathcal{B}}

\newcommand{\calV}{\mathcal{V}}

\newcommand{\dS}{{\mathfrak S}}

\newcommand{\bbS}{{\mathbb S}}
\newcommand{\bbP}{{\mathbb P}}
\newcommand{\bbR}{{\mathbb R}}

\newcommand{\Dch}{{\rm \bf D}}
\newcommand{\diag}{\operatorname{diag}}
\newcommand{\tr}{\operatorname{tr}}

\newcommand{\ch}{{\rm \bf Ch}}

\newcommand{\im}{\operatorname{im}}
\newcommand{\rank}{\operatorname{rank}}

\newtheorem{theorem}{Theorem}[section]
\newtheorem{prop}[theorem]{Proposition}
\newtheorem{corollary}[theorem]{Corollary}
\newtheorem{lm}[theorem]{Lemma}
\newtheorem{rem}[theorem]{Remark}
\newtheorem{example}[theorem]{Example}

\def\Proof{\noindent{\sl Proof.}\qquad}

\begin{document}

\title{Ranking patterns of unfolding models 
of codimension one}
\author{
Hidehiko Kamiya
\footnote
{
{\rm This work was partially supported by 
JSPS 
KAKENHI (19540131).} \,
{\it Faculty of Economics, Okayama University}
}
\\ 
Akimichi Takemura 
\footnote
{This research was supported by JST CREST. 
{\it Graduate School of Information Science and Technology,
University of Tokyo} 
} 
\\
Hiroaki Terao 
\footnote
{
{\rm This work was partially supported by 
JSPS KAKENHI (21340001).} \,
{\it Department of Mathematics, Hokkaido University}
}
} 
\date{August 2010}
\maketitle 

\begin{abstract}
We consider the problem of counting 
the number of possible sets of rankings 
(called ranking patterns)  generated by 
unfolding models of codimension one.  
We express the ranking patterns as slices of the 
braid arrangement and show that all braid slices, 
including those not associated with 
unfolding models, are in one-to-one correspondence with 
the chambers of an arrangement. 
By identifying those which are 
associated with unfolding models, 
we find the number of ranking patterns. 
We also give an upper bound for the number of 
ranking patterns  when the difference by a 
permutation  of objects is ignored. 
\end{abstract}


\section{Introduction} 
\label{sec:intro}


The unfolding model, 
also known as the ideal point model, 
is a model for preference rankings,    
and was introduced by 
Coombs \cite{coo-50} 
in 
psychometrics. 
Since then, this model has been widely used 
not only in psychometrics 
(De Soete, Feger and Klauer \cite{dfk}) 
but also in other 
fields such as marketing science 
(DeSarbo and Hoffman \cite{deh}, 
MacKay, Easley and Zinnes \cite{mez}).  
The same mathematical structure can also be found in 
voting theory (Hinich and Munger \cite{him}). 

In this paper, we consider the problem of 
counting 
the number of possible sets of 
rankings (called ranking patterns) 
generated by the unfolding model. 
We deal with the case where 
the restriction by dimension 
is weakest, 
and give the answer in terms of the number of chambers 
of a hyperplane arrangement. 

Suppose we have a set of $m$ objects 
labeled $1,2,\ldots,m$ 
and an individual who ranks these $m$ objects 
according to his/her preference.
In the unfolding model, it is assumed that 
the $m$ objects $1, 2, \ldots,m$ are represented by 
points $\mu_1, \mu_2, \ldots,\mu_m$ in 
the 
Euclidean space 
$\bbR^n$. 
Moreover, the individual is also 
represented by a point 
$y$ 
in 
the same $\bbR^n$. 
This 
$y$ is called the ideal point of the individual, 
and is identified with the individual. 
Then $\bbR^n$ containing both the objects and the 
individual is 
called the joint space in the 
psychometric 
literature. 
Now, according to the unfolding model, 
individual 
$y$ prefers object $i$ to object $j$ if and only if 
$y$ is closer to $\mu_i$ than to $\mu_j$ 
in the usual Euclidean distance, i.e., 
$\Vert y - \mu_i \Vert < \Vert y-\mu_j \Vert$.
So individual $y$ gives ranking 
$(i_1i_2 \cdots i_m)$, 
meaning that $i_1$ is the individual's best object, 
$i_2$ is his/her second best object, and so on, 
if and only if 
$y$ is closest to $\mu_{i_1}$, 
second closest to $\mu_{i_2}$, and so on.  

\begin{figure}[htbp]
 \begin{center}
\includegraphics*[width=.8\textwidth]{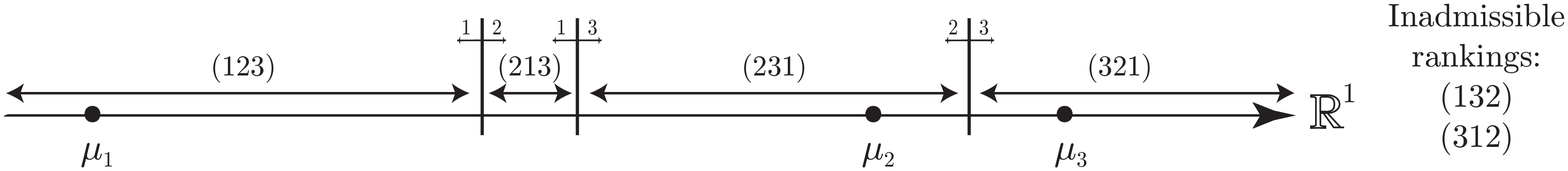}
\caption{$n=1, \ m=3$.}
\label{fig:n1m3}
 \end{center}
\end{figure}

In general, of course, we can think of $m!$ rankings 
among 
$m$ objects. 
But in the unfolding model, not all the $m!$ 
rankings are generated; 
there are admissible rankings and inadmissible rankings. 
That is, if there is a point $y$ in the joint space 
$\bbR^n$ 
which is closest to $\mu_{i_1}$, second closest to 
$\mu_{i_2}$, and so on, 
then the ranking $(i_1 i_2 \cdots i_m)$ is admissible. 
On the contrary, 
if there is no such point $y$, 
then 
$(i_1 i_2 \cdots i_m)$ 
is inadmissible. 
For the $m=3$ points $\mu_1, \mu_2, \mu_3$ on $\bbR^1$ 
in Figure \ref{fig:n1m3}, for instance, rankings 
$(123), (213), (231), (321)$ are admissible, while 
$(132), (312)$ are inadmissible.   
Then 
a natural question is: 
What is the number of admissible rankings 
for a given set of $m$ objects? 
For $n=1$, this number is obviously equal to 
$\binom{m}{2}+1$ as long as the midpoints of the objects 
are all distinct. 
But the 
question 
is not trivial for general $n\ge 2$.  
This problem has been solved, 
and the number is expressed in terms of  
the signless Stirling numbers of the first kind 
(Good and Tideman \cite{got}, 
Kamiya and Takemura \cite{kt-97, kt-05}, 
Zaslavsky \cite{zas-02}). 

Now, as we explained, 
for a given set of 
objects $\mu_1,\mu_2,\ldots,\mu_m$, 
we have admissible rankings. 
Let us call the set of all 
admissible rankings 
the ranking pattern of the unfolding model with 
$\mu_1,\mu_2,\ldots,\mu_m$. 
For the three objects $\mu_1,\mu_2,\mu_3$ in Figure \ref{fig:n1m3}, 
the ranking pattern is $\{ (123), (213), (231), (321)\}$.  
In general, 
if we change 
$\mu_1,\mu_2,\ldots,\mu_m$, 
we obtain a different ranking pattern. 
Our question is: 
How many ranking patterns are possible 
by taking different choices of 
$\mu_1,\mu_2,\ldots,\mu_m$? 

In the unidimensional case $n=1$,
determining the ranking pattern corresponds to 
determining the order of 
$m(m-1)/2$ 
midpoints of the 
objects on the real line $\bbR^1$. 
(See Lemma 2.3 and Theorem 2.4 of \cite{kott}.) 
For $m=3$ as in our example in Figure \ref{fig:n1m3}, 
there is only one possible order of midpoints 
if we restrict the order of objects as 
$\mu_1<\mu_2<\mu_3$. 
For general $m\ge 4$, however, there are many possible 
orders of midpoints, so counting this number is not easy. 
Thrall \cite{thr} gave an upper bound 
for 
the number of 
possible orders of midpoints, and thus the number of 
ranking patterns 
in the 
unidimensional case. 
He obtained his upper bound by considering a problem 
similar to that of counting the number of standard 
Young tableaux. 
Recently, Kamiya, Orlik, Takemura and Terao 
\cite{kott} found the exact number of 
ranking patterns of the unidimensional unfolding model. 
They showed that the exact number can be obtained 
by counting the number of chambers of 
an arrangement called 
the mid-hyperplane arrangement. 
(See also Stanley \cite{sta-06}.)
However, the problem of counting the number of 
ranking patterns is 
harder 
for 
general dimension. 

In the present paper, we consider the problem 
of counting the number of ranking patterns 
when the unfolding model is ``of codimension one,'' 
i.e., when $n=m-2$ 
so that the restriction by dimension is weakest. 
In this case, we show that there is a one-to-one 
correspondence between the set of ranking patterns 
and a subset of the set of chambers of an arrangement 
(a restriction of the ``all-subset arrangement''). 
By this one-to-one correspondence, we can obtain 
the number of ranking patterns.

Note that we distinguish the $m$ objects when we
count the number of ranking patterns.
We say that
two
ranking patterns 
are
equivalent 
when one is 
obtained from the other by a 
permutation 
of the objects. 
When we do not distinguish the objects,
we count 
the number of 
inequivalent ranking patterns. 
We give 
an upper bound for this number. 

The organization of this paper is as follows. 
In Section \ref{sec:um-as-bs}, 
we see 
that the ranking pattern of the unfolding model 
of codimension one can be 
obtained by slicing the braid arrangement 
by an affine hyperplane, 
although not all these slices
can be realized by unfolding models. 
In Section \ref{sec:rp-of-bs}, the set of 
braid slices is shown to be in one-to-one 
correspondence with 
the set of chambers of a restriction of 
the all-subset arrangement. 
Of 
these chambers, some correspond to 
braid slices 
realizable by unfolding models, 
and others correspond to unrealizable ones.
This distinction is made in 
Section \ref{sec:rbs}. 
Based on these results, we give the 
number of ranking patterns of unfolding models 
of codimension one 
in Section \ref{sec:n-of-rp-of-um}. 
In the final section (Section \ref{sec:irp}), 
we provide 
an upper 
bound for the number of inequivalent ranking patterns.

\section{The unfolding model as a braid slice}
\label{sec:um-as-bs}

In this section, we show that 
the ranking pattern of the unfolding model 
of codimension one can be 
obtained by slicing the braid arrangement 
by an affine hyperplane. 

Let $m$ be an integer with $m \ge 3$. 
Denote by $\bbP_m$  
the set of all 
permutations 
of $[m]:=\{1,\ldots,m\}$: 
$\bbP_m:=\{ (i_1\cdots i_m): 
(i_1 \cdots i_m) \text{ is a permutation of } [m]\}$. 

\begin{figure}[htbp]
 \begin{center}
\includegraphics*[width=.6\textwidth]{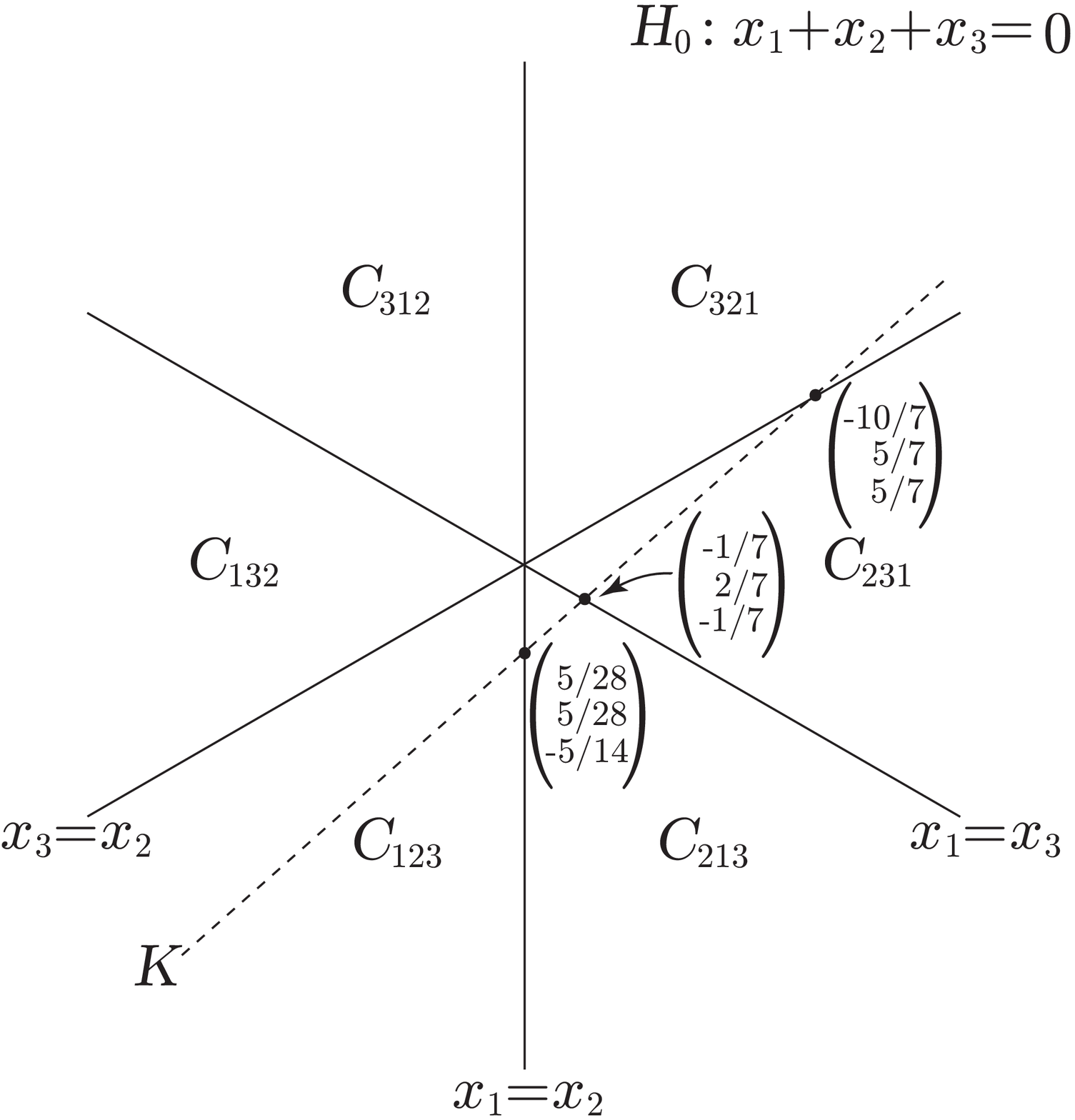}
\caption{$\calB_3^{H_0}$.}
\label{fig:B3}
 \end{center}
\end{figure}

Let 
\[
\calB_m := \{ H_{ij}: 1\le i<j \le m\}, \quad   
H_{ij} := \{ x=(x_1,\ldots,x_m)^T\in \bbR^m: 
x_i=x_j \}, 
\]
be the braid arrangement. 
Define 
\[
H_0:=\{ x=(x_1,\ldots,x_m)^T\in \bbR^m: 
x_1+\cdots +x_m=0\},  
\]
and set 
\[ 
C_{i_1\cdots i_m}
:=
\{ x=(x_1,\ldots,x_m)^T\in H_0: 
x_{i_1}>\cdots>x_{i_m}\}, 
\quad (i_1\cdots i_m)\in \bbP_m.  
\]
Note that 
$C_{i_1\cdots i_m}$ 
is a chamber of the arrangement 
$
\calB_m^{H_0}:=
\{ H\cap H_0: H\in \calB_m\}$ 
in $H_0$. 
Figure \ref{fig:B3} shows $\calB_3^{H_0}$. 

Now, for any $v\in 
\bbS^{m-2}:=\{ x \in H_0: \|x\|=1 \}$,  
let us define a hyperplane $K_{v}$ in $H_0$ by  
\[
K_{v}:=\{ x\in H_0: v^Tx=1\}. 
\]
We call 
\begin{equation}
\label{eq:RP(v)}
{\rm RP}(v):=\{ (i_1\cdots i_m)\in \bbP_m: 
K_{v}\cap C_{i_1\cdots i_m}\ne \emptyset \}, 
\quad v \in \bbS^{m-2}, 
\end{equation} 
the {\it ranking pattern of the braid slice} by 
$K_{v}$.

In general, for $m$ distinct points 
$\nu_{1}, \dots, \nu_{m} \in \bbR^{N} \ (m \ge N+1)$, 
let $\overline{\nu_{i}\nu_{j}}$ 
denote the one-simplex connecting two points
$\nu_{i}$ and $\nu_{j} \ (i < j)$.
Consider the following condition:


\begin{enumerate}
\item [(A)] 
The union of $N$ distinct one-simplices
$\overline{\nu_{i_{k} } \nu_{j_{k} } } $ 
$(i_{k} < j_{k}, \ k =1, \dots , N)$ contains no loop
if and only if the corresponding vectors 
$\nu_{i_{k} } - \nu_{j_{k} } $ 
$(k =1, \dots , N)$ 
are linearly independent.
\end{enumerate}


Recall, in general, 
that $N+1$ points $\tilde{\nu}_{1}, \dots, 
\tilde{\nu}_{N+1}  \in \bbR^{N}  $ 
are said to be
in general position if they are the vertices of an $N$-simplex,
in other words, the $N$ vectors 
$\tilde{\nu}_{1} - \tilde{\nu}_{2}, 
\tilde{\nu}_{2} - \tilde{\nu}_{3}, \dots, 
\tilde{\nu}_{N} - \tilde{\nu}_{N+1}$
are linearly independent.
It is not hard to see that condition (A)
implies that any $N+1$  points out of the 
$m$ points $\nu_1,\ldots,\nu_m$ 
are in general position.  
The converse,
however, is not true. 
For example,
$
\nu_{1}=(0, 0)^T, \   
\nu_{2}=(2, 0)^T, \  
\nu_{3}=(0, 1)^T
$ 
and
$\nu_{4}=(1, 1)^T$
do not satisfy condition (A)
because $\nu_{1} - \nu_{2}$ 
and
$\nu_{3} - \nu_{4} $ 
are linearly dependent, 
although any three of these $\nu_1,\ldots,\nu_4$ are 
in general position.  

Next, we 
move on to the ranking pattern of the 
unfolding model. 
Let $n\ge 1$ be a positive integer. 
By definition, $(i_1\cdots i_m)\in \bbP_m$ 
is admissible 
in the unfolding model with 
objects 
$\mu_1,\ldots,\mu_m \in \bbR^{n}$ 
iff there exists $y \in \bbR^{n}$ such that 
$\| y -\mu_{i_1}\|
< \cdots < \| y - \mu_{i_m} \|$.  
Let us call 
\begin{equation}
\label{eq:RP-UF}
{\rm RP}^{{\rm UF}}(\mu_1,\ldots,\mu_m)
:=
\{ (i_1\cdots i_m)\in \bbP_m: 
\| y -\mu_{i_1}\|< \cdots < \| y - \mu_{i_m} \| 
\text{ for some } y \in 
\bbR^n
\}
\end{equation}
the {\it ranking pattern of the unfolding model} with 
$\mu_1,\ldots,\mu_m\in 
\bbR^{n}$. 
Note that for every $c\in \bbR^n$, 
$\| y -\mu_{i_1}\|
< \cdots < \| y - \mu_{i_m} \|$ for some $y\in \bbR^n$ iff
$\| y -\mu_{i_1}-c\|
< \cdots < \| y - \mu_{i_m}-c \|$ for some $y\in \bbR^n$.
Hence the ranking pattern of the unfolding model is invariant with respect
to translations of $\mu_1,\dots,\mu_m$: 
\[
{\rm RP}^{{\rm UF}}(\mu_1,\ldots,\mu_m) = 
{\rm RP}^{{\rm UF}}(\mu_1+c,\ldots,\mu_m+c) ,  \quad c\in \bbR^n.
\]
Thus 
we can assume $\sum_{j=1}^m \mu_j=0_n$ without loss of generality, where
$0_{n}\in \bbR^{n}$ is the vector of zeros.    
As long as not all $\mu_1,\ldots,\mu_m$ are zero, 
we can also assume 
$
\sum_{j=1}^m ||\mu_j||^{2}/m = 1$ 
without loss of generality, 
because
the ranking pattern of the unfolding model is invariant with respect
to nonzero multiplications of $\mu_1,\dots,\mu_m$:
\[
{\rm RP}^{{\rm UF}}(\mu_1,\ldots,\mu_m) = 
{\rm RP}^{{\rm UF}}(a \mu_1,\ldots,a \mu_m) ,  
\quad a \in \bbR^{*}:=\bbR\setminus\{ 0\}.
\]
Therefore, we assume  
from now on 
that $\mu_1, \dots, \mu_m\in \bbR^n$ 
satisfy $\sum_{j=1}^m \mu_j=0_n$
and
$
\sum_{j=1}^m ||\mu_j||^{2}/m = 1$. 

Define an $m\times n$-matrix $W$ and 
an $m$-dimensional column vector $u$ by
\begin{eqnarray}
W
&=&
{\rm W}(\mu_1,\ldots,\mu_m)
=(w_1,\ldots,w_{n})
:=\begin{pmatrix}
\mu_1^T \\
\vdots \\  
\mu_m^T
\end{pmatrix}
\in {\rm Mat}_{m\times n}(\bbR), 
\label{eq:W}
\\ 
u
&=&
{\rm u}(\mu_1,\ldots,\mu_m):=
-\frac{1}{2}
\begin{pmatrix}
\|\mu_{1}\|^2 - 1
\\ 
\vdots \\ 
\|\mu_{m}\|^2 -1 
\end{pmatrix}
\in \bbR^m,  
\label{eq:u} 
\end{eqnarray}
where 
${\rm Mat}_{m\times n}(\bbR)$ stands for the 
set of $m\times n$-matrices with real entries. 
Consider an affine map
$\kappa : \bbR^{n} \longrightarrow \bbR^{m} $ defined
by
$\kappa(y) := Wy + u$ for $y\in\bbR^{n}$.
Let 
\[
K := \im \kappa=u + {\rm col \, } W, 
\]
where $\im \kappa:=\{ \kappa(y): y \in \bbR^n\}$ is 
the image of $\kappa$, and  
${\rm col \, } W$ is the column space of $W$. 
Note 
$1_m^TW=0_{n}^T$ and $1_m^Tu=0$  
(or $w_1,\ldots,w_{n}, u \in H_0$),  
where $1_m\in \bbR^m$ is the vector of ones.
Thus $K$ is an affine subspace of $H_{0} $. 
The condition defining 
${\rm RP}^{{\rm UF}}(\mu_1,\ldots,\mu_m)$ in 
\eqref{eq:RP-UF} can be expressed as follows: 
\begin{eqnarray}
&&
\| y -\mu_{i_1}\|< \cdots < \| y - \mu_{i_m} \| 
\text{ for some } y \in \bbR^{n} 
\label{eq:forsomey} \\ 
&& \qquad \iff 
\mu_{i_1}^T y
-\frac{1}{2}(\|\mu_{i_1}\|^2-1) 
>\cdots > 
\mu_{i_m}^T y
-\frac{1}{2}(\|\mu_{i_m}\|^2-1) 
\text{ for some } y \in \bbR^{n} 
\notag \\ 
&& \qquad \iff 
\kappa(y) \in C_{i_1\cdots i_m}
\text{ for some } y \in \bbR^{n} 
\notag \\ 
&& \qquad \iff 
K\cap C_{i_1\cdots i_m}
\ne \emptyset. 
\label{eq:KcapC} 
\end{eqnarray}
%
Condition \eqref{eq:KcapC} 
means that 
${\rm RP}^{{\rm UF}}(\mu_1,\ldots,\mu_m)$ 
can be obtained by slicing
the braid arrangement by an affine subspace. 
%

\begin{example}
\label{ex:n1m3}
Consider the case $n=1, \ m=3$, and take 
$\mu_1=-3/\sqrt{14/3}, \ \mu_2=1/\sqrt{14/3}, \ 
\mu_3=2/\sqrt{14/3}$. 
The objects in Figure \ref{fig:n1m3} were taken 
in this way. 
For these $\mu_1, \mu_2, \mu_3$, we have 
$W=(-3,1,2)^T/\sqrt{14/3}, \ u=(-13, 11, 2)^T/28$, 
so $K$ is the line 
defined by
\begin{equation}
\label{eq:line}
x_1=-\frac{13}{28}-3t, \ \ x_2=\frac{11}{28}+t, \ \ 
x_3=\frac{1}{14}+2t, \qquad t \in \bbR.  
\end{equation}
%
Find the line $K$ in Figure \ref{fig:B3} 
and compare Figure \ref{fig:B3} with Figure \ref{fig:n1m3}.
\end{example}

Consider the following two conditions on 
$\mu_1,\ldots,\mu_m\in \bbR^{n}
 \ (n\le m-2)$: 
\begin{enumerate}
\item [(A1)] 
The $m$ points $\mu_1,\ldots,\mu_m\in \bbR^{n} $
satisfy condition (A).
\item [(A2)] The $m$ points
$
\begin{pmatrix}
\mu_1
\\ 
\|\mu_{1} \|^{2}
\end{pmatrix},  
\ldots,
\begin{pmatrix}
\mu_m
\\ 
\|\mu_{m} \|^{2}
\end{pmatrix} 
\in \bbR^{n+1}$
satisfy condition (A).
\end{enumerate}
When 
$\mu_1,\ldots,\mu_m\in \bbR^{n}$ 
with $n\le m-2$ 
satisfy (A1) and (A2), we will say the unfolding model with $\mu_1,\ldots,\mu_m$ is 
(or $\mu_1,\ldots,\mu_m$ themselves 
are) {\it generic}. 
Note that (A1) and (A2) are translation invariant and 
nonzero multiplication invariant, 
i.e., $\mu_1, \dots, \mu_m$ are
generic iff 
$\mu_1+c, \dots, \mu_m+c$ 
are generic for any $c\in \bbR^n$
(or
$a \mu_1, \dots, a \mu_m$ 
are generic for any $a\in \bbR^{*}$).

\begin{rem}
When $n\ge m-1$, 
condition (A1) with the $N=n$ in (A) 
replaced by $m-1$ implies 
$\dim K=\rank W=m-1$ and thus $K=H_0$. 
In this case, $K\cap C_{i_1\cdots i_m}
=C_{i_1\cdots i_m}\ne \emptyset$ 
for all $(i_1\cdots i_m)\in \bbP_m$, 
and hence 
${\rm RP}^{{\rm UF}}(\mu_1,\ldots,\mu_m)$ 
is the whole $\bbP_m$: 
${\rm RP}^{{\rm UF}}(\mu_1,\ldots,\mu_m)
=\bbP_m$. 
\end{rem}
In the present paper, we will treat exclusively the 
case $n=m-2$.  

\begin{figure}[htbp]
 \begin{center}
\includegraphics*[width=.5\textwidth]{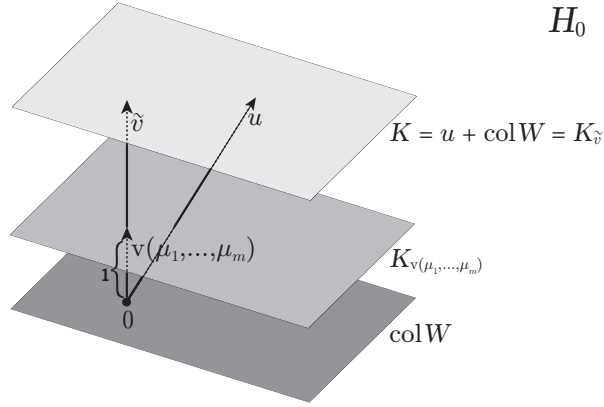}
\caption{$K_{{\rm v}(\mu_1,\ldots,\mu_m)}$.}
\label{fig:K_v}
 \end{center}
\end{figure}

Suppose 
$\mu_1,\ldots,\mu_m \in \bbR^{m-2}$ 
are generic. 
Let us define 
\begin{equation}
\label{eq:tilde-v}
\tilde{v}=\tilde{{\rm v}}(\mu_1,\ldots,\mu_m):=
u - {\rm proj}_{ {\rm col \, } W} (u), 
\end{equation}
where ${\rm proj}_{ {\rm col \, } W} $ stands for the
orthogonal projection on ${ {\rm col \, } W}$. 
Thanks to (A1), we have 
$\rank W=m-2$, 
so we can write $\tilde{v}$ as 
\begin{equation*}
\tilde{v}=
(I_m-W(W^TW)^{-1}W^T)u,
\end{equation*}
where $I_m$ denotes the identity matrix. 
Since the vector $u$ does not lie on ${ {\rm col \, } W} $
because of (A2), we have $\tilde{v}\ne 0_m$. 
Besides, we have $\dim K=m-2=\dim H_0-1$. 
These two facts imply that   
we can write $K=u + {\rm col \, } W$ in terms of 
$\tilde{v}$ as 
\[
K=K_{\tilde{v}}:=
\{ x\in H_0: \tilde{v}^Tx=\|\tilde{v}\|^2 \}. 
\]
Defining 
\begin{equation}
\label{eq:v}
{\rm v}(\mu_1,\ldots,\mu_m):=\frac{1}{\|\tilde{v}\|}
\tilde{v}, 
\end{equation}
we obtain the following equivalence: 
For $(i_1\cdots i_m)\in \bbP_m$, 
\begin{equation}
\label{eq:t-miffKcapC}
\| y -\mu_{i_1}\|< \cdots < \| y - \mu_{i_m} \| 
\text{ for some } y \in \bbR^{m-2}
\iff 
K_{{\rm v}(\mu_1,\ldots,\mu_m)}\cap C_{i_1\cdots i_m}\ne \emptyset, 
\end{equation}
where 
\[
K_{{\rm v}(\mu_1,\ldots,\mu_m)}
=\{ x\in H_0: {\rm v}(\mu_1,\ldots,\mu_m)^Tx=1\}. 
\] 
(See Figure \ref{fig:K_v}.) 

\begin{example}
\label{ex:Kv}
In the case of Example \ref{ex:n1m3}, we have 
$\tilde{v}=(5/98)(-1,5,-4)^T$ and hence 
$K_{{\rm v}(\mu_1,\mu_2,\mu_3)}=
\{ (x_1,x_2,x_3)^T\in H_0: -x_1+5x_2-4x_3=\sqrt{42}\}$, 
which is the dilate of line $K$ in \eqref{eq:line} by 
$\|\tilde{v}\|^{-1}=98/(5\sqrt{42})$.  
\end{example}

In the generic case with $n=m-2$, we have that 
$K$ 
is an affine hyperplane in $H_0$: 
$
\dim K
=\dim H_0-1,  
\ 0_m \notin 
K
$.  
We will say the unfolding model is 
{\it of codimension one} 
when $\mu_1,\ldots,\mu_m$ 
are generic 
with $n=m-2$.  

By \eqref{eq:RP(v)}, \eqref{eq:RP-UF} and 
\eqref{eq:t-miffKcapC},  
we obtain the following proposition. 

\begin{prop} 
\label{prop:rp-of-unfolding}
The ranking pattern of 
the 
unfolding model of codimension one 
with $\mu_1,\ldots,\mu_m\in \bbR^{m-2}$ 
is given by 
the ranking pattern 
of the braid slice 
by $K_{{\rm v}(\mu_1,\ldots,\mu_m)}$: 
\[
{\rm RP}^{{\rm UF}}(\mu_1,\ldots,\mu_m)=
{\rm RP}({\rm v}(\mu_1,\ldots,\mu_m)) 
\text{ for generic } 
\mu_1,\ldots,\mu_m \in \bbR^{m-2}.  
\]
\end{prop}

\section{Ranking patterns of braid slices}
\label{sec:rp-of-bs}

In this section, we show that the set of 
ranking patterns of 
braid slices by 
$K_{v} \ 
(v\in \bbS^{m-2})$
for ``generic'' $v$'s is in 
one-to-one correspondence with the set of chambers of 
an arrangement of hyperplanes in $H_0$. 
The discussions in this section are about 
braid slices, and 
the unfolding model does not 
concern us (except in a few places) in this section. 

We begin by defining 
an arrangement $\calA$ 
of hyperplanes in $\bbR^m$ by 
\[
\calA=\calA_m:=\{ H_I: I\subseteq [m], \ 
|I|\ge 1
\}, 
\quad 
H_I:=\{ x=(x_1,\ldots,x_m)^T\in \bbR^m: 
\sum_{i\in I}x_i=0\}. 
\]
Note that $H_0=H_{[m]}$. 
We will call $\calA$ the 
{\it all-subset arrangement}.  
Next we consider the restriction of $\calA$ to 
$H_0$: 
\[
\calA^0=\calA^0_{m}:=\calA_m^{H_0}=
\{ H^0_I: I \subset [m], \ 
1\le |I| \le m-1
\}, 
\quad 
H^0_{I}:=H_I\cap H_0. 
\] 
We notice that $H^0_{[m]\setminus I}=H^0_I$. 

Now define 
\[
\calV:=(H_0\setminus \bigcup \calA^0)\cap \bbS^{m-2}, 
\]
where $\bigcup \calA^0:=\bigcup_{H\in \calA^0}H$. 
Then we have the following basic lemma. 

\begin{lm} 
\label{lm:+--}
Take an arbitrary $v=(v_1,\ldots,v_m)^T \in \calV$. 
Then for $(i_1\cdots i_m)\in \bbP_m$,  
we have 
the equivalences below: 
\begin{eqnarray*}
&&
K_{v}\cap C_{i_1\cdots i_m}= \emptyset
\\ 
&& \qquad \iff v_{i_1}<0, \ v_{i_1}+v_{i_2}<0, \ldots, 
v_{i_1}+\cdots +v_{i_{m-1}}<0, 
\\ 
&&
K_{v}\cap C_{i_1\cdots i_m}\ne \emptyset  
\text{ is bounded}
\\
&& \qquad \iff v_{i_1}>0, \ v_{i_1}+v_{i_2}>0, \ldots, 
v_{i_1}+\cdots +v_{i_{m-1}}>0, 
\\  
&& 
K_{v}\cap C_{i_1\cdots i_m}
\text{ is unbounded}
\\ 
&& \qquad \iff \text{there exist } 
k,l\in [m-1] \ (k\ne l)  
\text{ such that } 
(v_{i_1}+\cdots +v_{i_k})(v_{i_1}+\cdots +v_{i_l})<0.
\end{eqnarray*}
\end{lm}

\Proof 
Without loss of generality, we may consider 
the case $(i_1\cdots i_m)=(1\cdots m)$. 
Let $c_1,\ldots,c_{m-1}\in H_0$ be 
defined by  
\begin{eqnarray*}
c_1&:=&(1,0,0,\ldots,0)^T-\frac{1}{m}1_m, \\ 
c_2&:=&(1,1,0,\ldots,0)^T-\frac{2}{m}1_m, \\ 
&& \vdots \\ 
c_{m-1}&:=&(1,1,\ldots,1,0)^T-\frac{m-1}{m}1_m. 
\end{eqnarray*}
Then $c_1,\ldots,c_{m-1}$ are linearly 
independent. 
Consider the pointed cone with apex $0_m$ and 
generators $c_1,\ldots,c_{m-1}$: 
\[
{\rm cone}\{ c_1,\ldots,c_{m-1}\}:=
\{ d_1c_1+\cdots +d_{m-1}c_{m-1}: 
d_1,\ldots,d_{m-1}\ge 0 \}, 
\]
which is a simplicial cone in $H_0$. 
Then $C_{1\cdots m}
=\{ (x_1,\ldots,x_m)^T\in H_0: 
x_1>\cdots >x_m \}$ is 
the relative interior of this cone: 
\[
C_{1\cdots m}=
{\rm relint}
({\rm cone}\{ c_1,\ldots,c_{m-1}\})
=\{ d_1c_1+\cdots +d_{m-1}c_{m-1}: 
d_1,\ldots,d_{m-1}>0 \}. 
\]

Suppose $v_1<0,\ v_1+v_2<0,\ldots,
v_1+\cdots +v_{m-1}<0$. 
This is equivalent to saying that 
$c_j^Tv<0$ for all $j\in [m-1]$, 
which in turn is equivalent to 
$K_{v}\cap 
{\rm relint}
({\rm cone}\{ c_1,\ldots,c_{m-1}\})
=\emptyset$. 

Suppose on the contrary that 
$v_1>0,\ v_1+v_2>0,\ldots,
v_1+\cdots +v_{m-1}>0$. 
Then $c_j^Tv>0, \ j\in [m-1]$, and 
hence we have 
\begin{equation}
\label{eq:relint=relint}
K_{v}\cap 
{\rm relint}
({\rm cone}\{ c_1,\ldots,c_{m-1}\})
=
{\rm relint}
\left(
{\rm conv}
\left\{ \frac{1}{c_1^Tv}c_1,\ldots,
\frac{1}{c_{m-1}^Tv}c_{m-1}\right\} 
\right),  
\end{equation}
where 
${\rm conv}\{ \ \}$ denotes the 
convex hull of the points in the braces.  
Noting that $c_1,\ldots,c_{m-1}$ are 
linearly independent, we can see that 
the right-hand side of 
\eqref{eq:relint=relint} is nonempty. 
Also, it is clearly bounded. 

Suppose instead that 
$v_1
+\cdots +v_k$ and 
$v_1
+\cdots +v_l$ have different signs 
for some $k$ and $l$. 
Then $c_k^Tv$ and $c_l^Tv$ have 
different signs. 
Hence, there exists $c\in{\rm relint}
({\rm cone}\{ c_1,\ldots,c_{m-1}\})
=C_{1\cdots m}
$ 
such that $c^Tv=0$. 
We have $v+dc\in K_{v}$ for any 
$d\in \bbR$; moreover, we can see  
$v+dc\in C_{1\cdots m}$ for all 
sufficiently large $d
>0$. 
Therefore, $K_{v}\cap C_{1\cdots m}$ is 
an unbounded set.  

Since there are no other cases than 
the three above 
for the signs of $\sum_{j=1}^sv_j \ 
(s\in [m-1])$ for $v=(v_1,\ldots,v_m)^T\in \calV$, 
the preceding arguments 
suffice 
to prove the three equivalences in the lemma. 
\qed

\

By \eqref{eq:RP(v)} and Lemma \ref{lm:+--}, 
it is easily seen that 
$|\bbP_m \setminus {\rm RP}(v)|=(m-1)!$ for 
any $v\in \calV$. 
When ${\rm RP}(v)$ can be realized 
by the unfolding model, 
this 
follows also 
from the general result 
on the cardinality of a ranking pattern 
of the unfolding model 
(Good and Tideman \cite{got}, 
Kamiya and Takemura \cite{kt-97, kt-05}, 
Zaslavsky \cite{zas-02}).

Let $\ch(\calA^0)$ stand for the set of chambers of 
$\calA^0$. 
Then we can write $\calV$ as 
\[
\calV
=\bigsqcup_{\tilde{D}\in \ch(\calA^0)}
(\tilde{D}\cap \bbS^{m-2}) 
=\bigsqcup_{D\in \Dch(\calA^0)}D
\quad \text{ (disjoint union)}, 
\]
where 
\[
\Dch(\calA^0):=\{ D=\tilde{D}\cap\bbS^{m-2}: 
\tilde{D}\in \ch(\calA^0)\}
\] 
is in 
one-to-one correspondence with 
$\ch(\calA^0)$. 
Using Lemma \ref{lm:+--}, we can prove the 
following proposition. 

\begin{prop}
\label{prop:chA-RP}
There is a one-to-one correspondence between 
$\Dch(\calA^0)$ 
and $\{{\rm RP}(v): v\in \calV\}$ 
given by the bijection 
\begin{equation}
\label{eq:ChtoRP}
\Dch(\calA^0)\ni D 
\longmapsto 
{\rm RP}(v), \quad 
v\in 
D. 
\end{equation}
\end{prop}

\Proof 
It is clear that the map \eqref{eq:ChtoRP} is 
well-defined and surjective. 
We will show that it is injective. 
Suppose 
$D$ and $D'$ ($D, D' \in \Dch(\calA^0)$) 
are different.  
Take arbitrary 
$v\in D$ 
and 
$v'\in D'$. 
Then there exists 
$I\subset [m], \ 1\le | I | \le m-1$, such that  
$\sum_{i\in I}v_i$ and $\sum_{i\in I}v'_i$ have 
different signs. 
Without loss of generality, we may assume 
$\sum_{i\in I}v_i<0$ and $\sum_{i\in I}v'_i>0$.  
Define $I_-=\{ i \in I: v_i<0\} \ne \emptyset, \ 
I_+
=\{ i \in I: v_i>0\}, \ 
\bar{I}_-
=\{ i\in [m]\setminus I: v_i<0 \}, \ 
\bar{I}_+
=\{ i\in [m]\setminus I: v_i>0 \} \ne \emptyset$. 
Take an arbitrary $(i_1\cdots i_m)\in \bbP_m$ such 
that $\{ i_1,\ldots,i_{|I_-|}\}=I_-, \ 
\{ i_{|I_-|+1},\ldots,i_{ |I| }\}=I_+, \ 
\{ i_{|I|+1},\ldots,i_{|I|+|\bar{I}_-|}\}=\bar{I}_-, \ 
\{ i_{|I|+|\bar{I}_-|+1},\ldots,i_{[m]}\}=\bar{I}_+$. 
Then $v=(v_i,\ldots,v_m)^T$ satisfies $v_{i_1}<0, \ 
v_{i_1}+v_{i_2}<0,\ldots,v_{i_1}
+\cdots +v_{i_{m-1}}<0$. 
Thus we have $(i_1\cdots i_m)\notin {\rm RP}(v)$ by 
Lemma \ref{lm:+--}. 
On the other hand, this is not the case with 
$v'=(v'_1,\ldots,v'_m)^T$ 
because $v'_{i_1}+\cdots+v'_{i_{|I|}}=
\sum_{i\in I}v'_i>0$, 
and we have 
$(i_1\cdots i_m)\in {\rm RP}(v')$ by 
Lemma \ref{lm:+--}. 
Therefore, we obtain ${\rm RP}(v)\ne {\rm RP}(v')$. 
\qed

\


Proposition \ref{prop:chA-RP} implies that 
the ranking patterns 
${\rm RP}(v), \ v \in \calV=
\bigsqcup_{D\in \Dch(\calA^0)}D$,   
are the same on a common $D$ and 
different on different $D$'s.   
So we can write ${\rm RP}(v)$ with $v\in D$ as 
${\rm RP}_D$: 
\[ 
{\rm RP}_D:={\rm RP}(v), \quad v\in D 
\in \Dch(\calA^0),   
\] 
and we have ${\rm RP}_D\ne {\rm RP}_{D'}$ for 
$D\ne D'$. 
%
%

We will say the braid slice by $K_{v} \ 
(v\in \bbS^{m-2})$ is {\it generic} when 
$v\in \calV$.   
It can be checked that if 
$\mu_1,\ldots,\mu_m\in \bbR^{m-2}$ are generic, 
${\rm v}(\mu_1,\ldots,\mu_m)$ defined in 
\eqref{eq:v} satisfies ${\rm v}(\mu_1,\ldots,\mu_m) 
\in \calV$, i.e., 
the braid slice by $K_{{\rm v}(\mu_1,\ldots,\mu_m)}$ is 
generic.

\section{Realizable braid slices}
\label{sec:rbs}

By Proposition \ref{prop:rp-of-unfolding}, we know 
that the ranking pattern of any 
unfolding model 
of codimension one 
can be obtained as the ranking pattern of 
a 
generic braid slice. 
However, not all 
ranking patterns 
of generic braid slices, 
${\rm RP}_{D}, \ D\in \Dch(\calA^0)$, 
can be realized 
as 
ranking patterns of 
unfolding models 
of codimension one. 
In this section, we establish
conditions on 
$D\in \Dch(\calA^0)$ 
which guarantee that 
${\rm RP}_D$ 
can be realized by 
an 
unfolding model 
of codimension one. 

Let $\calV_2$ be the set of all 
$v
\in \calV$ having 
at least two positive entries 
and at the same time at least two negative entries: 
\begin{eqnarray*}
\calV_2 &:=& \{ v=(v_1,\ldots,v_m)^T \in \calV: 
v_i,v_j>0 \text{ and } v_k,v_l<0 \\ 
&& \qquad \qquad \qquad \qquad \qquad \ \ \text{ for some } 
i,j,k,l\in [m] \ (i\ne j, \ k\ne l)
\}. 
\end{eqnarray*}
Then put $\calV_1:=\calV\setminus \calV_2$. 
We see that $\calV_1$ is the set of all 
$v=(v_1,\ldots,v_m)^T\in \calV$ having 
exactly one positive entry or exactly one negative 
entry. 
Note that 
for any $D\in \Dch(\calA^0)$, we have either 
$D\subset \calV_2$ or $D\subset \calV_1$.  



It is helpful to consider 
$D\in \Dch(\calA^0)$ and $-D=\{ -v: v\in D\} 
\in \Dch(\calA^0)$ in a pair. 
Obviously, $D\subset \calV_i$ implies 
$-D\subset \calV_i$ for each $i=1,2$. 

\begin{theorem}
\label{th:both+exactly-one}
For any $D\in \Dch(\calA^0)$, we have the 
following.  
\begin{enumerate}
\item 
Suppose $D\subset \calV_2$. 
Then each of ${\rm RP}_D$ and ${\rm RP}_{-D}$ can be 
realized as the ranking pattern of 
an unfolding model of codimension one, 
i.e., there exist generic 
$\mu_1,\ldots, \mu_m\in \bbR^{m-2}$ and 
$\mu'_1,\ldots, \mu'_m\in \bbR^{m-2}$ such that 
\[
{\rm RP}_D
={\rm RP}^{{\rm UF}}(\mu_1,\ldots, \mu_m), \quad 
{\rm RP}_{-D}
={\rm RP}^{{\rm UF}}(\mu'_1,\ldots, \mu'_m). 
\]

\item Suppose $D\subset \calV_1$. 
Then  exactly one 
of ${\rm RP}_D$ and ${\rm RP}_{-D}$ can be 
realized as the ranking pattern of 
an 
unfolding model 
of codimension one. 
In fact, ${\rm RP}_{\varepsilon D}$ can be realized and 
${\rm RP}_{-\varepsilon D}$ cannot be realized, where 
$\varepsilon =\pm 1$
is such that 
$\varepsilon v$ for any $v\in D$ 
has exactly one positive entry. 
\end{enumerate}
\end{theorem}


The proof of Theorem \ref{th:both+exactly-one} is 
based on the following two lemmas. 
For $v=(v_1,\ldots,v_m)^T\in\bbR^m$, 
let  
$\diag(v):=\diag(v_1,\ldots,v_m) \in {\rm Mat}_{m\times m}(\bbR)$
stand for the diagonal matrix with diagonal 
entries 
$v_1, \dots, v_m$.

\begin{lm}
\label{lm:WdW}
Suppose $v\in \bbR^m$ and 
$W\in {\rm Mat}_{m\times (m-2)}(\bbR)$ 
satisfy  
\[
v\ne 0_m, \quad 1_m^Tv=0, \quad 
1_m^TW=v^TW=0_{m-2}^T. 
\]
\begin{enumerate}
\item 
If $v$ has at least two positive entries 
as well as 
at least two negative entries, 
then $W^T\diag(v)W$ is indefinite (i.e., has 
at least one positive eigenvalue and at least one 
negative eigenvalue).  
\item If $v$ has exactly one positive 
(resp. negative) entry, 
then $W^T\diag(v)W$ is non-positive (resp. 
non-negative) definite.  
If in addition $v$ has at least two 
negative (resp. positive) entries, then 
$W^T\diag(v)W$ has at least one 
negative (resp. positive) 
eigenvalue, and hence $\tr\{ W^T\diag(v)W\}$ is 
negative (resp. positive). 
\end{enumerate}
\end{lm}

\Proof 
We can assume without loss of generality that $v
=(v_1,\ldots,v_m)^T$ is of unit length: 
$\|v\|^2=\sum_{i=1}^mv_i^2=1$.  
Define $C:=(1_m, 
v-2^{-1}\sum_{i=1}^mv_i^3 1_m, 
W-1_m(v_1^2,\ldots,v_m^2)W)\in 
{\rm GL}(m,\bbR)$.  
Then, by direct calculations, we can see that 
\begin{equation}
\label{eq:CdvC}
C^T\diag(v)C=
\begin{pmatrix}
0 & 1 & 0_{m-2}^T \\ 
1 & 0 & 0_{m-2}^T \\ 
\ \ \quad 
0_{m-2} & 
\ \ \quad 
0_{m-2} & 
W^T\diag(v)W
\end{pmatrix}. 
\end{equation}
Equation \eqref{eq:CdvC} implies that 
the number of positive (resp. negative) eigenvalues of 
$W^T\diag(v)W$ plus one 
is equal to 
the number of positive (resp. negative) eigenvalues of 
$C^T\diag(v)C$, 
which in turn is equal to 
the number of positive (resp. negative) entries of $v$ 
by Sylvester's law of inertia. 
\qed


\begin{lm}
\label{lm:ABB}
Suppose an $m'\times m'$ real symmetric matrix $A$ 
is indefinite.  
Then we have 
$\{ \tr(BA B^T): B\in {\rm GL}(m',\bbR) \}=\bbR$.  
\end{lm}

\Proof 
Let $\lambda_1,\ldots,\lambda_{m'}$ be the eigenvalues 
of $A$ with $\lambda_1>0$ and $\lambda_2<0$, 
and write $\Lambda=\diag(\lambda_1,\ldots,
\lambda_{m'})$. 
Then 
\begin{eqnarray*}
\{ \tr (B A B^T): B\in {\rm GL}(m',\bbR) \}
&=& 
\{ \tr(B\Lambda B^T): B\in {\rm GL}(m',\bbR) \} \\ 
&=&
\{ \lambda_1 \|b_1\|^2+\lambda_2 \|b_2\|^2
+\lambda_3 \|b_3\|^2
+\cdots + \lambda_{m'}\|b_{m'}\|^2: \\ 
&& \qquad \qquad \qquad \qquad \qquad 
(b_1,\ldots,b_{m'})\in {\rm GL}(m',\bbR)
\}. 
\end{eqnarray*}
For given $b_3,\ldots,b_{m'}$, we can take 
$b_1, b_2$ with arbitrary positive 
lengths. 
\qed 
 
\ 

\noindent {\sl 
Proof of Theorem \ref{th:both+exactly-one}.}  
\quad 
Take an arbitrary 
$v
\in D$.  
Let $\{ w_1,\ldots,w_{m-2}\}$ be a basis of 
$H_0\cap ({\rm span}\{v\})^{\perp}=
({\rm span}\{ 1_m, v\})^{\perp}$, 
and take $\mu_1,\ldots,\mu_m\in \bbR^{m-2}$ as  
\[
\begin{pmatrix}
\mu_1^T \\ 
\vdots \\ 
\mu_m^T
\end{pmatrix}
:=(w_1,\ldots,w_{m-2}). 
\]  
Note that 
$\sum_{j=1}^m\mu_j=0_{m-2}$. 
Moreover, we can take $w_1,\ldots,w_{m-2}$ 
so that $\sum_{j=1}^m\| \mu_j \|^2/m=1$.  
For such $\mu_1,\dots,\mu_m$, 
let us consider 
${\rm W}(\mu_1,\dots,\mu_m)$ 
and 
${\rm u}(\mu_1,\dots,\mu_m)$ 
defined in 
\eqref{eq:W} and \eqref{eq:u}: 
\[
{\rm W}(\mu_1,\ldots,\mu_m)
=
\begin{pmatrix}
\mu_1^T \\ 
\vdots \\ 
\mu_m^T
\end{pmatrix}
\in {\rm Mat}_{m\times (m-2)}(\bbR), 
\quad 
{\rm u}(\mu_1,\ldots,\mu_m)=
-\frac{1}{2}
\begin{pmatrix}
\|\mu_1\|^2-1 \\ 
\vdots \\ 
\|\mu_m\|^2-1 
\end{pmatrix}
\in \bbR^m. 
\]
We note here that ${\rm u}(\mu_1,\ldots,\mu_m)^Tv$ 
can be written as 
\begin{equation}
\label{eq:umuv}
{\rm u}(\mu_1,\ldots,\mu_m)^Tv
=-\frac{1}{2}\tr\{ 
{\rm W}(\mu_1,\ldots,\mu_m)^T
\diag(v)
{\rm W}(\mu_1,\ldots,\mu_m)\}. 
\end{equation}
Moreover, using the fact that 
$v\notin \bigcup \calA^0$, we can check 
that $\mu_1,\ldots,\mu_m$ satisfy (A1). 

We first prove Part 1. 
Suppose $D\subset \calV_2$. 
Then, since $v\in D \subset \calV_2$, 
we have by 
Part 1 of 
Lemma \ref{lm:WdW} that  
the symmetric matrix 
${\rm W}(\mu_1,\ldots,\mu_m)^T\diag(v)
{\rm W}(\mu_1,\ldots,\mu_m)$ is indefinite. 
So Lemma \ref{lm:ABB} implies that 
there exist $B_1\in {\rm GL}(m-2,\bbR)$ 
and $B_2\in {\rm GL}(m-2,\bbR)$ 
such that 
\begin{gather*}
\tr\{ B_1{\rm W}(\mu_1,\ldots,\mu_m)^T\diag(v)
{\rm W}(\mu_1,\ldots,\mu_m)B_1^T\}>0, \\ 
\tr\{ B_2{\rm W}(\mu_1,\ldots,\mu_m)^T\diag(v)
{\rm W}(\mu_1,\ldots,\mu_m)B_2^T\}<0. 
\end{gather*}
Together with 
\begin{eqnarray*}
{\rm u}(B_k\mu_1,\ldots,B_k\mu_m)^Tv
&=&
-\frac{1}{2}\tr\{ 
{\rm W}(B_k\mu_1,\ldots,B_k\mu_m)^T
\diag(v){\rm W}(B_k\mu_1,\ldots,B_k\mu_m)\} \\ 
&=& 
-\frac{1}{2}\tr\{ 
B_k
{\rm W}(\mu_1,\ldots,\mu_m)^T
\diag(v)
{\rm W}(\mu_1,\ldots,\mu_m)
B_k^T\},  
\quad k=1,2, 
\end{eqnarray*}
these inequalities imply 
\begin{equation}
\label{eq:uBmv}
{\rm u}(B_1\mu_1,\ldots,B_1\mu_m)^Tv<0, \quad   
{\rm u}(B_2\mu_1,\ldots,B_2\mu_m)^Tv>0. 
\end{equation}
We observe that the column space of 
${\rm W}(B_k\mu_1,\ldots,B_k\mu_m)=
{\rm W}(\mu_1,\ldots,\mu_m)B_k^T$ 
is equal to 
that of 
${\rm W}(\mu_1,\ldots,\mu_m)$.  
This fact and 
${\rm u}(B_k\mu_1,\ldots,B_k\mu_m)^Tv \ne 0$ 
yield 
\[
{\rm v}(B_k\mu_1,\ldots,B_k\mu_m)
={\rm sign}
\{ {\rm u}(B_k\mu_1,\ldots,B_k\mu_m)^Tv\}v, 
\quad 
k=1,2
\]
(see \eqref{eq:tilde-v}, \eqref{eq:v} and 
Figure \ref{fig:K_v}). 
By \eqref{eq:uBmv}, we obtain 
\[
{\rm v}(B_1\mu_1,\ldots,B_1\mu_m)=-v, 
\quad  
{\rm v}(B_2\mu_1,\ldots,B_2\mu_m)=v. 
\]
Now, since $\mu_1,\ldots,\mu_m$ satisfy (A1), 
clearly so do 
$B_k\mu_1,\ldots,B_k\mu_m$ 
for $k=1,2$.  
From this fact and 
${\rm u}(B_k\mu_1,\ldots,B_k\mu_m)^Tv\ne 0$, 
we can check that 
$B_k\mu_1,\ldots,B_k\mu_m$ also 
satisfy (A2) ($k=1,2$). 
Now that $B_k\mu_1,\ldots,B_k\mu_m$ are generic 
($k=1,2$), 
Proposition \ref{prop:rp-of-unfolding} yields  
${\rm RP}(-v)=
{\rm RP}({\rm v}(B_1\mu_1,\ldots,B_1\mu_m))
={\rm RP}^{{\rm UF}}(B_1\mu_1,\ldots,B_1\mu_m)$ 
and  
${\rm RP}(v)=
{\rm RP}({\rm v}(B_2\mu_1,\ldots,B_2\mu_m))
={\rm RP}^{{\rm UF}}(B_2\mu_1,\ldots,B_2\mu_m)$.  
Thus, we have proved that each of 
${\rm RP}(v)$ 
and ${\rm RP}(-v)$ 
is realized by 
an 
unfolding model 
of codimension one, 
where $v\in D$ and $-v\in -D$. 
This completes the proof of Part 1. 

Next we prove Part 2. 
Suppose $D\subset \calV_1$. 
Then the fact that 
$v\in \calV_1$ 
together with 
Part 2 of Lemma \ref{lm:WdW} 
and 
equation 
\eqref{eq:umuv} 
implies that  
${\rm u}(\mu_1,\ldots,\mu_m)^Tv
\ne 0$. 
Hence we have 
${\rm v}(\mu_1,\ldots,\mu_m)
=\varepsilon v, \ 
\varepsilon={\rm sign}
\{ {\rm u}(\mu_1,\ldots,\mu_m)^Tv\}$.  
Also, 
from ${\rm u}(\mu_1,\ldots,\mu_m)^Tv
\ne 0$ and the fact that 
$\mu_1,\ldots,\mu_m$ satisfy (A1),  
it follows that 
$\mu_1,\ldots,\mu_m$ satisfy (A2) 
as well.  
Thus we obtain 
${\rm RP}(\varepsilon v)=
{\rm RP}({\rm v}(\mu_1,\ldots,\mu_m))
={\rm RP}^{{\rm UF}}
(\mu_1,\ldots,\mu_m)$. 
This proves that at least one of 
${\rm RP}(v)$ and ${\rm RP}(-v)$ can 
be realized by an unfolding 
model of codimension one. 
It remains to show that not both 
${\rm RP}(v)$ and ${\rm RP}(-v)$ can 
be realized by unfolding models of codimension one. 
Suppose on the contrary that 
both ${\rm RP}(v)$ and ${\rm RP}(-v)$ 
were realized. 
Without loss of generality, assume that
$v_{i_0}<0, \ v_i>0 \ (i\ne i_0)$ for 
some $i_0\in [m]$, where $v_i \ 
(1\le i\le m)$ are the entries of $v$.   
But by taking
$y$ in \eqref{eq:forsomey} sufficiently 
close to $\mu_{i_0}$, 
we see that ${\rm RP}(v)$ with such a 
$v$ cannot  be realized by an 
unfolding model of codimension one, 
because ${\rm RP}(v)=\bbP_m\setminus 
\{ (i_0i_1\cdots i_{m-1}): 
(i_1\cdots i_{m-1}) \text{ is a 
permutation of } [m]\setminus \{ i_0\}\}$. 
This is a contradiction. 
\qed

\section{The number of ranking patterns of 
unfolding models}
\label{sec:n-of-rp-of-um}

Based on the results in 
Sections \ref{sec:um-as-bs}, 
\ref{sec:rp-of-bs} and \ref{sec:rbs}, 
we find, in this section, 
the number of ranking patterns of 
unfolding models of codimension one. 

For $i \in [m]$, let us define $\calV_1(i,+)
\subset \calV_1$ by 
\[
\calV_1(i,+):=
\{ v=(v_1,\ldots,v_m)^T\in \calV_1: 
v_{i}>0, \ v_{j}<0 
\text{ for all } 
j \in [m]\setminus \{i\}\}. 
\]

\begin{lm}
\label{lm:V=D} 
For any $i \in [m]$, 
we have $\calV_1(i,+)=D_i$ for some 
$D_i\in \Dch(\calA^0)$. 
\end{lm}

\Proof 
Obviously, $\calV_1(i,+)$ is a union of some 
chambers 
$D\in \Dch(\calA^0)$. 
So it suffices to show the 
following: 
For any $I\subset [m]$ with $1\le |I|\le m-1$, we 
have $\calV_1(i,+)\subset (H_I^0)^{+}
\cap \bbS^{m-2}$  
or $\calV_1(i,+)\subset (H_I^0)^{-}\cap \bbS^{m-2}$, 
where 
$(H_I^0)^{+}:=\{ x=(x_1,\ldots,x_m)^T
\in H_0: \sum_{j\in I}x_j>0\}$ and 
$(H_I^0)^{-}:=H_0\setminus(H_I^0\cup (H_I^0)^+)$.  
If $i\notin I$, any $v=(v_1,\ldots,v_m)^T
\in \calV_1(i,+)$ satisfies $\sum_{j\in I}v_j<0$, 
and thus we have 
$\calV_1(i,+)\subset (H_I^0)^{-}
\cap \bbS^{m-2}$. 
If $i\in I$, on the other hand, 
$v=(v_1,\ldots,v_m)^T
\in \calV_1(i,+)$ implies $\sum_{j\in I}v_j
=-\sum_{j\in [m]\setminus I}v_j>0$, 
so we obtain $\calV_1(i,+)\subset 
(H_I^0)^{+}\cap \bbS^{m-2}$. 
\qed 

\ 

We can write $\calV_1$ as 
\[
\calV_1=
D_1\sqcup (-D_1)
\sqcup \cdots \sqcup 
D_m\sqcup (-D_m),  
\]
where 
$-D_i = \{ -v: v\in D_i\}\in \Dch(\calA^0)$ 
for $i\in [m]$.  
Notice $-D_i = \calV_1(i,-)$ with   
\[
\calV_1(i,-):=
\{ v=(v_1,\ldots,v_m)^T\in \calV_1: 
v_{i}<0, \ v_{j}>0 
\text{ for all } 
j \in [m]\setminus \{i\}\}  
\] 
for $i \in [m]$.  

Now, consider the mapping 
\begin{eqnarray*}
{\rm v}: 
\{ 
(\mu_1,\ldots,\mu_m): 
\mu_1,\ldots,\mu_m \in \bbR^{m-2} \text{ are generic} 
\} 
&\longrightarrow& 
\calV, \\ 
(\mu_1,\ldots,\mu_m) 
&\longmapsto& {\rm v}(\mu_1,\ldots,\mu_m). 
\end{eqnarray*}
From the proof of Theorem \ref{th:both+exactly-one}, 
we can see that the image 
$
\im {\rm v}=
\{ 
{\rm v}(\mu_1,\ldots,\mu_m): 
\mu_1,\ldots, \mu_m \in \bbR^{m-2} 
\text{ are generic} 
\}
$
of 
${\rm v}$ 
is given by 
\begin{eqnarray}
\im {\rm v} 
&=&
\bigsqcup_{D\in \Dch(\calA^0), \ 
D \ne 
-D_i 
\ (i\in [m])} D 
\notag \\ 
&=&
\calV 
\setminus 
((-D_1) \sqcup \cdots \sqcup (-D_m))
=
\calV_2 
\sqcup D_1 \sqcup \cdots \sqcup D_m. 
\label{eq:im-v}  
\end{eqnarray}

We are in a position to state the 
main result of this section. 
Denote by $q(m)$ 
the number of ranking patterns of 
unfolding models of codimension one: 
\[
q(m):=|\{{\rm RP}^{{\rm UF}}(\mu_1,\ldots,\mu_m): 
\text{ generic } \mu_1,\ldots,\mu_m \in \bbR^{m-2}\}|.  
\]

\begin{theorem}
\label{th:q(m)}
It holds that 
$
q(m)=
|\ch(\calA^0)|-m, \ \calA^0=\calA^0_m. 
$
\end{theorem}

\Proof 
By Propositions \ref{prop:rp-of-unfolding} and 
\ref{prop:chA-RP} 
and equation \eqref{eq:im-v}, we have 
\[
q(m)=  
|\{
{\rm RP}(v): 
v \in 
\calV 
\setminus ((-D_1)\sqcup \cdots \sqcup(-D_m))
\}| 
=
|\Dch(\calA^0)|-m. 
\]
\qed 

\

We have calculated specific values of $q(m)$ for 
$m\le 8$ 
in the following way. 

The number of chambers 
$|\ch(\calA_m^0)|$ can be 
obtained by finding the characteristic polynomial 
$\chi(\calA^0_m, t)$ of $\calA_m^0$ 
(Orlik and Terao \cite[Definition 2.52]{ort}): 
$|\ch(\calA_m^0)|
=(-1)^{m-1}\chi(\calA^0_m, -1)$ 
(Zaslavsky \cite[Theorem A]{zas-75}, 
Orlik and Terao \cite[Theorem 2.68]{ort}).  
Moreover, when finding $\chi(\calA^0_m, t)$, 
we can use the property 
$L(\calA_m^0) \simeq L(\calA_{m-1})$ of the 
all-subset arrangement,  
where $L( \, \cdot \, )$ denotes the 
intersection poset of an arrangement 
(Orlik and Terao \cite[Definition 2.1]{ort}).  
The characteristic polynomials 
$\chi(\calA^0_m, t)$ and 
the numbers of chambers 
$|\ch(\calA^0_m)|$ of $\calA^0_m$ for $m\le 8$ 
are given in the following lemma. 

\begin{lm}
\label{lm:chi(3-8)}
For $m \le 8$, 
$\chi(\calA^0_m, t)$ and 
$|\ch(\calA^0_m)|$ are given by 
\begin{eqnarray*}
\chi(\calA^0_3, t)
&=& t^2-3t+2=(t-1)(t-2), \quad 
|\ch(\calA^0_3)|= 6; \\ 
\chi(\calA^0_4, t)
&=& t^3-7t^2+15t-9=(t-1)(t-3)^2, \quad 
|\ch(\calA^0_4)|= 32; \\ 
\chi(\calA^0_5, t)
&=& 
t^4-15t^3+80t^2-170t+104 
=(t-1)(t-4)(t^2-10t+26), \\ 
&& |\ch(\calA^0_5)|
= 370; \\ 
\chi(\calA^0_6, t)
&=& t^5-31t^4+375t^3-2130t^2+5270t-3485 \\ 
&=& (t-1)(t^4-30t^3+345t^2-1785t+3485), \\ 
&& |\ch(\calA^0_6)|
= 11292; \\ 
\chi(\calA^0_7, t)
&=& t^6-63t^5+1652t^4-22435t^3+159460t^2-510524t+371909 \\ 
&=& (t-1)(t^5-62t^4+1590t^3-20845t^2+138615t-371909), \\ 
&& |\ch(\calA^0_7)|
= 1066044; \\ 
\chi(\calA^0_8, t)
&=& t^7-127t^6+7035t^5-215439t^4+3831835t^3 \\ 
&& \qquad \qquad \qquad 
-37769977t^2+169824305t-135677633 \\ 
&=& (t-1)(t^6-126t^5+6909t^4-208530t^3 \\ 
&& \qquad \qquad \qquad \qquad 
+3623305t^2-34146672t+135677633
), \\ 
&& |\ch(\calA^0_8)|
= 347326352. 
\end{eqnarray*} 
\end{lm}
%
%
%

We obtained $\chi(\calA^0_3, t)$ and 
$\chi(\calA^0_4, t)$ by direct calculations. 
For $\chi(\calA^0_5, t)$, we used the method of 
deletion and restriction 
(Orlik and Terao \cite[Theorem 2.56]{ort}).     
Furthermore, we calculated 
$\chi(\calA^0_6, t), \ \chi(\calA^0_7, t)$ and 
$\chi(\calA^0_8, t)$ 
by the finite field method 
(Athanasiadis \cite{ath-thesis, ath-96}, 
Stanley \cite[Lecture 5]{sta-07}, 
Crapo and Rota \cite{crr}, 
Kamiya, Takemura and Terao \cite{ktt-08, ktt-root, 
ktt-nc}). 


\begin{rem} 
We can consider $\calA_m^0$ also for $m=2$, and we have 
$\chi(\calA^0_2, t)=t-1$ and $|\ch(\calA^0_2)|=2$. 
The arrangement $\calA^0_m \ (m\ge 2)$ 
also appears in 
thermal field theory 
(Evans \cite{eva-92, eva-95}, van Eijck \cite{vei}). 
The numbers $|\ch(\calA^0_m)| \ (m=2,\ldots,8)$: 
\[
2, \ 6, \ 32, \ 370, \ 11292, \ 1066044, \ 347326352
\]
are listed in \cite[Table 1]{eva-95} and 
\cite[Table 2.1]{vei} 
as the numbers of 
regions 
of the analytic continuations 
of ITF (imaginary-time formalism) 
Green functions,  
although the characteristic polynomials 
$\chi(\ch(\calA^0_m), t) \ (m \le 8)$ 
are not obtained there.   
\end{rem}

From Theorem \ref{th:q(m)} and the values of 
$
|\ch(\calA^0_m)| \ (3 \le m \le 8)$ 
in Lemma \ref{lm:chi(3-8)}, 
we can obtain 
$q(m) \ (3 \le m \le 8)$:  

\begin{corollary}
The numbers $r(m)$ 
of ranking patterns of unfolding models of codimension 
one for $m\le 8$ are given by 
\begin{gather*}
q(3)
=3, \ 
q(4)
=28, \ 
q(5)
=365, \\ 
q(6)
=11286, \ 
q(7)
=1066037, \ 
q(8)
=347326344. 
\end{gather*}
\end{corollary}

%
%

\section{Inequivalent ranking patterns} 
\label{sec:irp}

In this section, we define 
equivalence of ranking patterns, 
and give an upper bound for the number of 
inequivalent ranking patterns of 
unfolding models of codimension one. 
For $m\le 6$,  we will see that this upper bound 
is actually the exact number. 

\subsection{The number of inequivalent 
ranking patterns of unfolding models}

Let $\dS_m$ be the symmetric group 
on $m$ letters, 
consisting of all bijections 
$\sigma:[m] \longrightarrow
[m]$.  
Let us say that ranking patterns 
${\rm RP}_D$ and 
${\rm RP}_{D'} \ (D, D'\in \Dch(\calA^0))$ 
of generic braid slices are {\it equivalent} 
iff 
\[
{\rm RP}_D=\sigma {\rm RP}_{D'} 
\text{ for some }  
\sigma \in \dS_m, 
\]
where 
\[
\sigma {\rm RP}_{D'}:=
\{ (\sigma(i_1)\cdots \sigma(i_m)): 
(i_1 \cdots i_m)\in {\rm RP}_{D'}\}. 
\]
We say ${\rm RP}_D$ and ${\rm RP}_{D'}$ are 
{\it inequivalent} iff they are not equivalent. 
We want to compute 
the number of 
inequivalent ranking patterns of generic 
braid slices 
that can be realized by 
unfolding models of codimension one. 

Consider the action of $\dS_m$ on 
$\calV$ 
defined by 
\[
\dS_m \times \calV 
\ni (\sigma, v) \longmapsto 
\sigma v:=(v_{\sigma^{-1}(1)},\ldots,v_{\sigma^{-1}(m)})^T
\in \calV, 
\]
where $v=(v_1,\ldots,v_m)^T$. 
This induces the action of $\dS_m$ on $\Dch(\calA^0)$: 
\begin{eqnarray}
\label{eq:SDtoD}
\dS_m \times \Dch(\calA^0)\ni (\sigma, D) \longmapsto 
\sigma D:=\{ \sigma v: v \in D\} 
\in \Dch(\calA^0). 
\end{eqnarray}
%
We can check 
\[
{\rm RP}_{\sigma D}=\sigma {\rm RP}_{D},  
\quad D \in \Dch(\calA^0), 
\ \sigma \in \dS_m.
\] 
Thus, ${\rm RP}_D$ and 
${\rm RP}_{D'}$ 
are equivalent 
iff $D$ and $D'$ are on the same orbit 
under 
action \eqref{eq:SDtoD}. 
Therefore, the number of 
inequivalent ranking patterns of generic 
braid slices is equal to 
the number of orbits 
$\dS_mD:=\{ \sigma D: \sigma \in \dS_m \}, 
\ D\in \Dch(\calA^0)$,  
i.e., the cardinality of the orbit space 
$\Dch(\calA^0)/\dS_m:=
\{ \dS_mD: D\in \Dch(\calA^0)\}$ 
under action \eqref{eq:SDtoD}. 

For each orbit 
$\dS_mD \in \Dch(\calA^0)/\dS_m$, 
either all or none of its elements 
$D'\in \dS_mD$ correspond to ranking patterns 
${\rm RP}_{D'}$ 
realizable 
by unfolding models of codimension one. 
Among 
the orbits in $\Dch(\calA^0)/\dS_m$, 
exactly one orbit, 
$\dS_m\calV_1(m, -)=\dS_m(-D_m)$, 
consists of elements (chambers) that 
correspond to ranking patterns 
not realizable 
by unfolding models of codimension one, 
${\rm RP}_{-D_1},\ldots,{\rm RP}_{-D_m}:$ 
\[
\dS_m\calV_1(m, -)
=\{ -D_1,\ldots,-D_m\}
\] 
(see \eqref{eq:im-v}). 
Therefore, the number of inequivalent 
ranking patterns 
${\rm RP}_{D} \ (D \in \Dch(\calA^0))$ 
realizable 
by unfolding models of codimension one is 
$|\Dch(\calA^0)/\dS_m|-1$. 


For ranking patterns of unfolding 
models of codimension one,  
we say 
${\rm RP}^{{\rm UF}}(\mu_1,\ldots,\mu_m)$ 
and 
${\rm RP}^{{\rm UF}}(\mu'_1,\ldots,\mu'_m)$ 
are {\it equivalent} (resp. {\it inequivalent}) 
iff they are equivalent (resp. inequivalent) 
when regarded as ranking patterns of 
generic braid slices. 
So 
${\rm RP}^{{\rm UF}}(\mu_1,\ldots,\mu_m) 
={\rm RP}({\rm v}(\mu_1,\ldots,\mu_m))$  
and 
${\rm RP}^{{\rm UF}}(\mu'_1,\ldots,\mu'_m)
={\rm RP}({\rm v}(\mu'_1,\ldots,\mu'_m))$ 
are equivalent iff 
$D\ni {\rm v}(\mu_1,\ldots,\mu_m)$ and 
$D'\ni {\rm v}(\mu'_1,\ldots,\mu'_m)$ 
are on the same orbit 
under action \eqref{eq:SDtoD}. 
Of course, if 
$(\mu_1,\ldots,\mu_m)=
(\mu'_{\sigma^{-1}(1)},\ldots,
\mu'_{\sigma^{-1}(m)})$ 
for some $\sigma\in \dS_m$, then 
${\rm RP}^{{\rm UF}}(\mu_1,\ldots,\mu_m)$ 
and 
${\rm RP}^{{\rm UF}}(\mu'_1,\ldots,\mu'_m)$ 
are equivalent 
(because of 
${\rm v}(\mu'_{\sigma^{-1}(1)},\ldots,
\mu'_{\sigma^{-1}(m)})
=\sigma {\rm v}(\mu'_1,\ldots,\mu'_m)$ 
yielding 
${\rm RP}^{{\rm UF}}(\mu_1,\ldots,\mu_m)
=\sigma {\rm RP}^{{\rm UF}}(\mu'_1,\ldots,\mu'_m)$), 
but not vice versa.  

From the arguments above, we obtain the 
following proposition. 

\begin{prop}
\label{prop:D/S-1}
The number of inequivalent ranking patterns 
of unfolding models of codimension one is 
equal to 
$|\Dch(\calA^0)/\dS_m|-1$. 
\end{prop}

Finding $|\Dch(\calA^0)/\dS_m|$ is not always easy.  
We will give an upper bound for 
the number, $|\Dch(\calA^0)/\dS_m|-1$,  
of inequivalent ranking patterns 
of unfolding models of codimension one. 

We have 
$D \not\subset 
\bigcup \calB_m^{H_0}
$ for any $D\in \Dch(\calA^0)$.  
Thus,  
to each orbit $\dS_mD \in \Dch(\calA^0)/\dS_m$, 
there belongs  
a chamber $\sigma D \in \dS_mD 
\ (\sigma \in \dS_m)$  
that 
intersects  
$C_{1\cdots m}
$.  
Hence, 
the set 
\begin{equation*}
\Dch^{1\cdots m}(\calA^0):=
\{ D\in \Dch(\calA^0): 
D\cap C_{1\cdots m}\ne \emptyset\}
\end{equation*} 
always includes 
a cross section 
(i.e., a complete set of representatives 
of the orbits) 
under action \eqref{eq:SDtoD}. 
Therefore, 
an 
upper bound for $|\Dch(\calA^0)/\dS_m|$ is given by 
the cardinality of 
$\Dch^{1\cdots m}(\calA^0)$: 
\begin{equation}
\label{eq:upper}
|\Dch(\calA^0)/\dS_m| 
\le 
|\Dch^{1\cdots m}(\calA^0)| 
=\frac{|\ch(\calA^0\cup \calB_m^{H_0})|}{m!}.  
\end{equation}
If, in particular, $\Dch^{1\cdots m}(\calA^0)$ 
is a cross section, 
then the inequality 
in \eqref{eq:upper} is actually 
an 
equality. 

Note that we can write 
$\Dch^{1\cdots m}(\calA^0)$ 
as 
\begin{equation}
\label{eq:{}={}cup{}}
\Dch^{1\cdots m}(\calA^0)
= \{ D_1, -D_m\} 
\cup 
\Dch_2^{1\cdots m}(\calA^0), 
\end{equation}
where 
\[
\Dch_2^{1\cdots m}(\calA^0):=
\{ D\in \Dch(\calA^0): 
D\subset \calV_2, \  
D\cap C_{1\cdots m}\ne \emptyset\}. 
\] 
Note, moreover, that for each 
$D \in 
\Dch_2^{1\cdots m}(\calA^0)$, 
three 
chambers $D_1,-D_m$ and $D$ 
are all on different orbits.  
Thus   
$\Dch^{1\cdots m}(\calA^0)$ 
is a cross section 
if and only if all elements of 
$
\Dch_2^{1\cdots m}(\calA^0)
$ 
are on different orbits. 
Define $\rho \in \dS_m$ by  
\begin{equation}
\label{eq:rho}
\rho(i):=m+1-i, \quad i\in [m]. 
\end{equation} 
Then, for any $D \in 
\Dch_2^{1\cdots m}(\calA^0)$, 
we have that 
$D$ and 
$-\rho D
=\{ (-v_{m},\ldots,-v_{1})^T: 
(v_1,\ldots,v_m)^T\in D\}
\in 
\Dch_2^{1\cdots m}(\calA^0)$ 
are on different orbits.  This can be seen as follows.
Without loss of generality, suppose $v_1 > \dots > v_m$
and $v_1+v_m > 0$.  Then  $v_1+v_j > 0$ for all $j=2,\dots,m$. But then
there is no $i$ such that
$-v_i - v_j > 0$ for all $j\neq i$.
However, the fact that $D$ and $-\rho D$ are on different orbits 
does not exclude the 
possibility 
of some $D$ and $D'$ 
($D, D'\in \Dch_2^{1\cdots m}(\calA^0), 
\ D\ne D'$) being on the same orbit.    


For $\Dch_2^{1\cdots m}(\calA^0)$,  
we may find 
$\Dch_2^{\overline{1\cdots m}}(\calA^0)$ 
instead: 
$\Dch_2^{1\cdots m}(\calA^0)
=\Dch_2^{\overline{1\cdots m}}(\calA^0)$, 
where 
$\Dch_2^{\overline{1\cdots m}}(\calA^0):=
\{ D\in \Dch(\calA^0): 
D \subset \calV_2, \ 
D\cap \bar{C}_{1\cdots m}\ne \emptyset\}$ 
with 
$\bar{C}_{1\cdots m}:=
\{ (x_1,\ldots,x_m)^T\in H_0: x_1\ge \cdots \ge x_m\}$. 

By Proposition \ref{prop:D/S-1}, 
\eqref{eq:upper} and \eqref{eq:{}={}cup{}}, 
we obtain an upper bound 
for 
the number of inequivalent 
ranking patterns of unfolding models of 
codimension one.  

\begin{corollary}
The number of inequivalent ranking patterns 
of unfolding models of codimension one 
cannot exceed 
\begin{equation}
\label{eq:1+||}
1+|
\Dch_2^{1\cdots m}(\calA^0)
| 
= \frac{|\ch(\calA^0\cup \calB_m^{H_0})|}{m!}-1. 
\end{equation} 
Moreover, 
if all elements of 
$
\Dch_2^{1\cdots m}(\calA^0)
$ 
are on different orbits 
under action \eqref{eq:SDtoD},  
then \eqref{eq:1+||} gives the exact number 
of inequivalent ranking patterns 
of unfolding models of codimension one. 
\end{corollary}
%
%
%
%


\subsection{Inequivalent ranking patterns 
for $m\le 6$} 

In this subsection, we investigate inequivalent 
ranking patterns of unfolding models of codimension one 
for $m\le 6$. 

We know 
\begin{eqnarray*}
{\rm RP}_{D_1}
&=& 
\bbP_m\setminus 
\{(i_1\cdots i_{m-1} \, 1): (i_1 \cdots i_{m-1}) 
\text{ is a permutation of } 
\{ 2,\ldots,m \}\}, \\ 
{\rm RP}_{-D_m}
&=&{\rm RP}_{-\rho D_1} 
=
\bbP_m\setminus 
\{(m \, i_1\cdots i_{m-1}): (i_1\cdots i_{m-1}) 
\text{ is a permutation of } 
[m-1]\}
\notag 
\end{eqnarray*}
by Lemma \ref{lm:+--}.

\subsubsection{Case $m=3$} 

When $m=3$, we have $\calV_2=\emptyset$ and 
$\Dch(\calA^0_3)=\{ D_1, D_2, D_3, -D_1, -D_2, -D_3 \}$. 
Accordingly, the set of all ranking patterns of 
unfolding models of codimension one is 
$\{ {\rm RP}_{D_1}, {\rm RP}_{D_2}, {\rm RP}_{D_3}\}$. 
Since $D_1=\tau_2 D_2=\tau_3 D_3$ 
($\tau_2\in \dS_3$ is the transposition of $1$ and $2$, 
and 
$\tau_3\in \dS_3$ is the transposition of $1$ and $3$), 
the number of inequivalent ranking patterns of 
unfolding models of codimension one 
is $1$. 
We have 
${\rm RP}_{D_1}=\bbP_3\setminus 
\{(231), (321)\}$. 

Let us consider 
\eqref{eq:upper} 
in this case. 
We have 
$\Dch_2^{123}(\calA_3^0)=\emptyset$, and 
$
\Dch^{123}(\calA_3^0)=
\{ D_1, -D_3\}\subset \Dch (\calA_3^0)$ 
is a cross section under the action of $\dS_3$ on 
$\Dch (\calA^0_3)$:  
\[
\Dch (\calA^0_3)=\dS_3 D_1 \sqcup 
\dS_3(-D_3)
=\{ D_1, D_2, D_3\} \sqcup 
\{ -D_1, -D_2, -D_3 \}.  
\]
Thus, the inequality in \eqref{eq:upper} 
is actually an 
equality 
in this case: 
$|\Dch(\calA^0_3)/\dS_3|
=
|
\Dch^{123}(\calA_3^0)
|
=2$. 
The number 
$|
\Dch^{123}(\calA_3^0)
|
=2$ 
can also be confirmed by 
$\chi(\calA^0_3\cup\calB_3^{H_0}, t)
=t^2-6t+5=(t-1)(t-5)$ 
yielding 
$|\ch(\calA^0_3\cup \calB_3^{H_0})|/(3!)
=(-1)^{3-1}\chi(\calA^0_3\cup\calB_3^{H_0}, -1)/(3!)
=12/(3!)=2$. 



\subsubsection{Case $m=4$} 

When $m=4$, we have 
\begin{equation}
\label{eq:V2capC1234}
\calV_2 
\cap \bar{C}_{1\cdots 4} 
=\bigsqcup_{
D\subset \calV_2, \ 
D \in \Dch(\calA^0_4)
}
(D\cap \bar{C}_{1\cdots 4})
=R_{4}\sqcup (-\rho 
R_{4}), 
\end{equation}
where 
\begin{equation*}
R_{4}:=
\{ (v_1,v_2,v_3,v_4)^T\in \bbS^2: 
v_1\ge v_2 >0> v_3 \ge v_4, \ v_2>-v_3\},  
\end{equation*} 
$\rho \in \dS_4$ is defined in \eqref{eq:rho} 
and 
\begin{eqnarray*}
-\rho R_4
&:=&\{ -(v_{\rho^{-1}(1)},\ldots,
v_{\rho^{-1}(4)})^T: 
(v_1,\ldots,v_4)^T\in R_4\} \\ 
&=&\{ (v_1,v_2,v_3,v_4)^T\in \bbS^2: 
v_1\ge v_2 >0> v_3 \ge v_4, \ v_2<-v_3\}.   
\end{eqnarray*}  

Now, there is only one 
$D\subset \calV_2 \ 
(D \in \Dch(\calA^0_4))$ such that 
$\emptyset\ne 
D\cap \bar{C}_{1\cdots 4}\subseteq R_4$. 
Such a $D$ is the chamber $E 
\ (\ni 
(v_1,\ldots,v_4)^T)$ 
determined by 
$v_{\{ 1\}}, v_{\{ 2\}}, v_{\{ 1, 2\}}, v_{\{ 1, 3\}}, v_{\{ 2, 3\}}
>0$ 
and 
$v_{\{ 3\}}, v_{\{ 4\}}
<0$,   
and we have 
$\emptyset \ne E \cap \bar{C}_{1\cdots 4}
=R_4$.  
Here, we are writing $v_I:=\sum_{i\in I}v_i$. 
As for $-\rho R_4$ in \eqref{eq:V2capC1234}, 
$D=-\rho E$ 
is the only 
$D\subset \calV_2 \ 
(D \in \Dch(\calA^0_4))$ such that 
$\emptyset\ne 
D\cap \bar{C}_{1\cdots 4}\subseteq -\rho R_4:  
\emptyset \ne (-\rho E)\cap \bar{C}_{1\cdots 4}
=-\rho R_4$.  

From the preceding arguments, we obtain 
\[
\calV_2 
\cap \bar{C}_{1\cdots 4} 
=
(E\cap \bar{C}_{1\cdots 4})
\sqcup 
((-\rho E)\cap \bar{C}_{1\cdots 4}), 
\] 
and hence 
$\Dch_2^{1\cdots 4}(\calA^0)
=\Dch_2^{\overline{1\cdots 4}}(\calA^0)
=\{ E, -\rho E\}$. 
Thus, we get  
\begin{equation*}
\Dch^{1\cdots 4}(\calA_4^0)
=\{ D_1, -\rho D_1, E, -\rho E\}
\end{equation*}
by \eqref{eq:{}={}cup{}}. 
We know that 
$
E$ 
and $-\rho E$ are 
on different orbits.   
Therefore, 
$\Dch^{1\cdots 4}(\calA_4^0)$ 
is a cross section 
under the action of $\dS_4$ on 
$\Dch(\calA_4^0)$, 
and we have $|\Dch(\calA^0_4)/\dS_4|
=|\Dch^{1\cdots 4}(\calA_4^0)|
=2(1+1)=4$.  
So the number of inequivalent ranking patterns 
of unfolding models of codimension one 
is $4-1=3=1+2\cdot 1$.  
(In passing, we can confirm the number 
$|\Dch^{1\cdots 4}(\calA_4^0)|=4$ by 
$\chi(\calA^0_4\cup\calB_4^{H_0}, t)
=t^3-13t^2+47t-35=(t-1)(t-5)(t-7)$ 
giving 
$|\ch(\calA^0_4\cup \calB_4^{H_0})|/(4!)
=(-1)^{4-1}\chi(\calA^0_4\cup\calB_4^{H_0}, -1)/(4!)
=96/(4!)=4$.) 
 
\begin{figure}[htbp]
 \begin{center}
\includegraphics*[width=.8\textwidth]{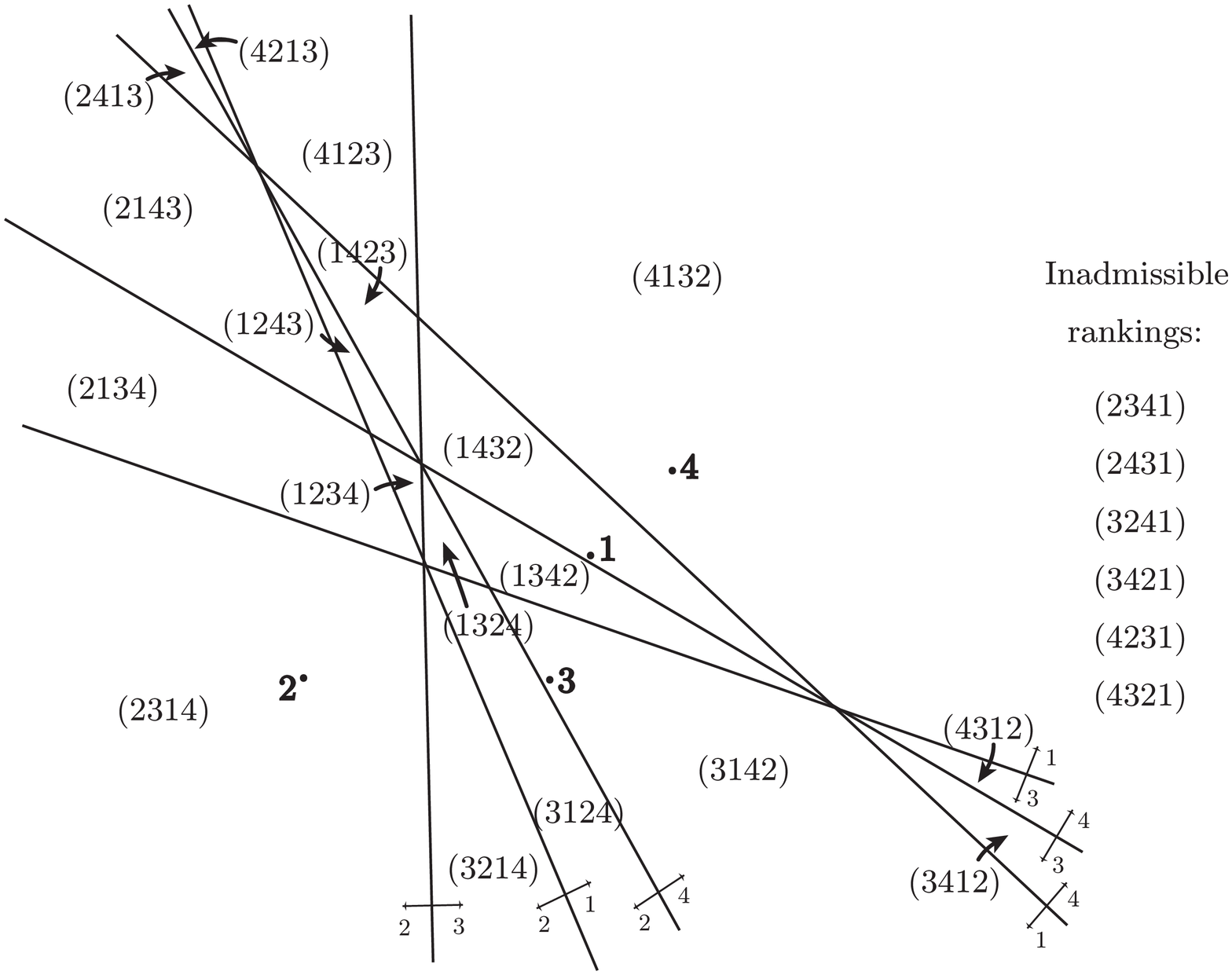}
\caption{${\rm RP}_{D_1}$.}
\label{fig:D1}
 \end{center}
\end{figure}

\begin{figure}[htbp]
 \begin{center}
\includegraphics*[width=.8\textwidth]{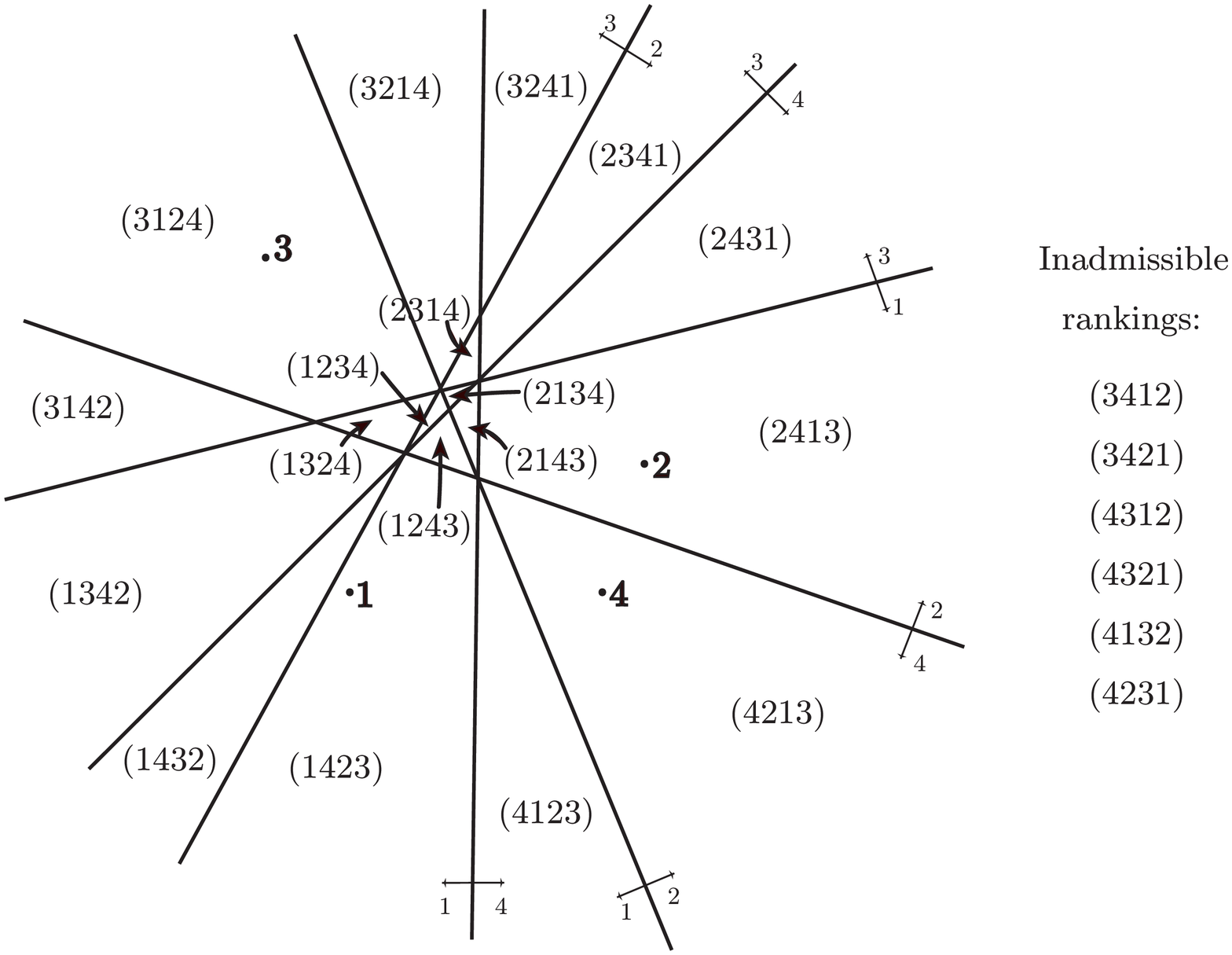}
\caption{${\rm RP}_{E}$.}
\label{fig:E}
 \end{center}
\end{figure}

\begin{figure}[htbp]
 \begin{center}
\includegraphics*[width=.8\textwidth]{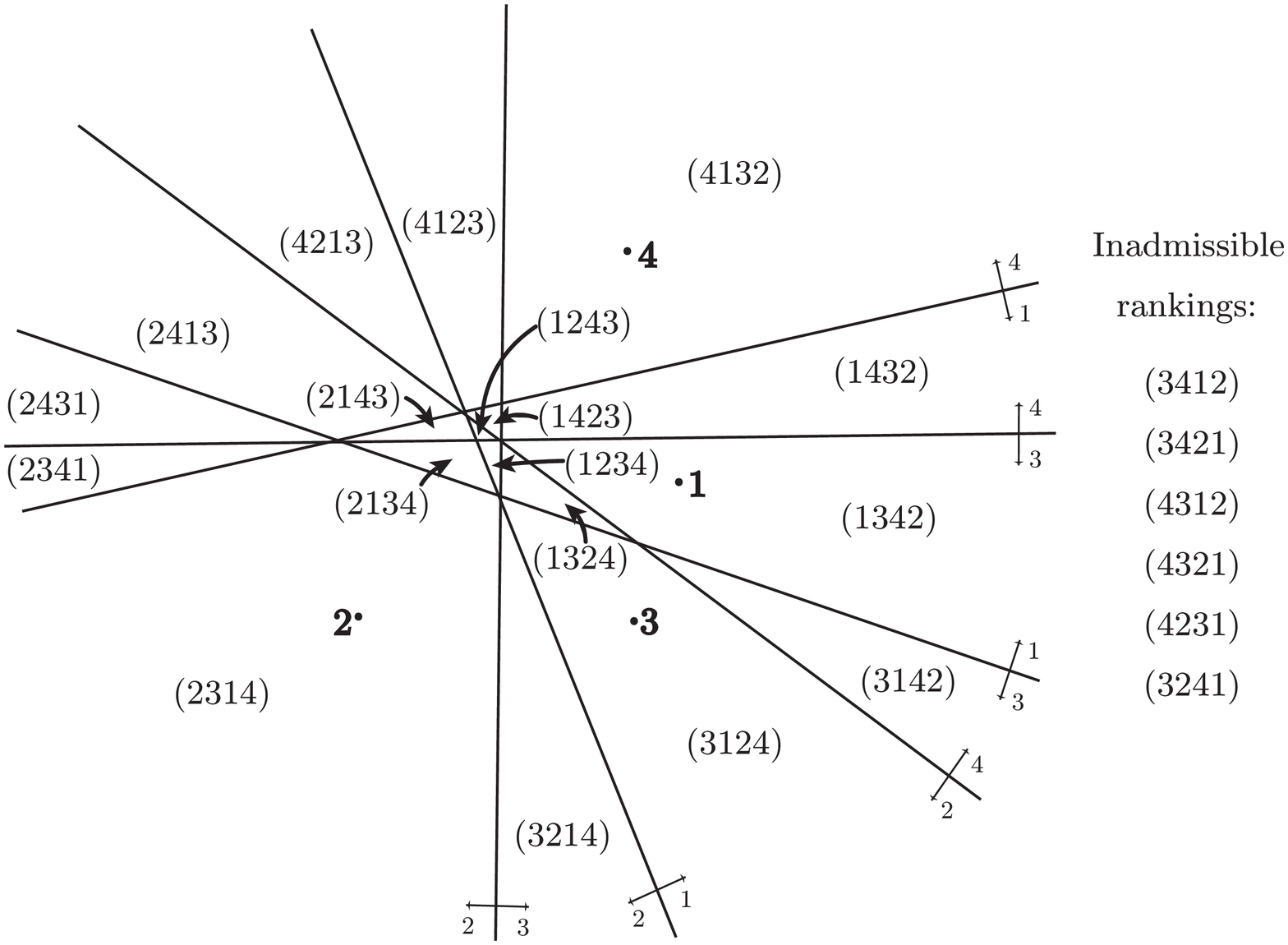}
\caption{${\rm RP}_{-\rho E}$.}
\label{fig:-rE}
 \end{center}
\end{figure}

The chambers $D_1, E, -\rho E$ 
correspond to ranking patterns 
that can be realized by unfolding 
models of codimension one, 
${\rm RP}_{D_1}, {\rm RP}_{E}, 
{\rm RP}_{-\rho E}$, 
while the chamber 
$-\rho D_1=-D_4=\calV_1(4,-)
$ corresponds to  
${\rm RP}_{-D_4}$, which cannot be realized.
From $E$, we can take 
$v=(1/2, 1/2, -1/4, -3/4)^T/(3\sqrt{2}/4)
\in E$. 
Thus by Lemma \ref{lm:+--}, we can see 
\begin{eqnarray*}
{\rm RP}_{D_1}
&=& 
\bbP_m \setminus 
\{ (2341), (2431), (3241), (3421), 
(4231), (4321)\}, \\ 
{\rm RP}_{E}
&=& 
\bbP_m \setminus 
\{ (3412), (3421), (4312), (4321), 
(4132), (4231)\}, \\ 
{\rm RP}_{-\rho E}
&=& 
\bbP_m \setminus 
\{ (3412), (3421), (4312), (4321), 
(4231), (3241)\}. 
\end{eqnarray*}
These three ranking patterns,  
realized as the ranking patterns of 
unfolding models of codimension one,  
are displayed 
in Figures \ref{fig:D1}, \ref{fig:E} and \ref{fig:-rE}. 
(For simplicity, $\mu_1,\mu_2, \mu_3, \mu_4$ 
are written as $1,2,3,4$ in the figures.) 

\subsubsection{Case $m=5$} 


When $m=5$, we have 
\begin{equation}
\label{eq:V2capC1-5}
\calV_2\cap \bar{C}_{1\cdots 5}
=\bigsqcup_{
D\subset \calV_2, \ 
D \in \Dch(\calA^0_5)
}
(D\cap \bar{C}_{1\cdots 5})
=R_{5}\sqcup(-\rho R_{5}), 
\end{equation} 
where  
\[
R_5:=\{ (v_1,\ldots,v_5)^T\in \bbS^3: 
v_1\ge v_2 \ge v_3 > 0 > v_4 \ge v_5\}. 
\]

There are five 
chambers 
$D\subset \calV_2 \ 
(D \in \Dch(\calA^0_5))$ such that 
$\emptyset\ne 
D\cap \bar{C}_{1\cdots 5}\subseteq R_5$. 
Let $E_1,\ldots,E_5$ be those 
five 
chambers. 
All  
$D\cap \bar{C}_{1\cdots 5}$ for 
$D=E_1,\ldots,E_5$ 
are listed in Table \ref{table:R5}. 
The first column gives the defining inequalities of 
each $D\cap \bar{C}_{1\cdots 5}$ 
(besides $v_1 \ge v_2 \ge v_3 > 0 > v_4 \ge v_5$); 
the second column exhibits an unnormalized 
representative point of  
each $D\cap \bar{C}_{1\cdots 5}$; and 
the last column contains an upper bound 
for 
$|\dS_5D|$ for each $D$. 
For example, 
the first row of Table \ref{table:R5} 
corresponds to the chamber $E_1 
\ni v=
(v_1,\ldots,v_5)^T
$ 
determined by 
positive 
$v_{\{ 1\}}, v_{\{ 2\}}, v_{\{ 3\}}, 
v_{\{ 1,2\}}, v_{\{ 1,3\}}, 
v_{\{ 1,4\}},
v_{\{ 2,3,4\}}, 
v_{\{ 2,3\}}, 
v_{\{ 2,4\}}, 
v_{\{ 1,3,4\}}, 
v_{\{ 3,4\}}, 
v_{\{ 1,2,4\}} 
$ 
and 
negative 
$v_{\{ 4\}}, v_{\{ 5\}}, 
v_{\{ 4,5\}}
$. 
For this 
$E_1$, 
we can take  
$v=(1/3, 1/3, 1/3, -1/6, -5/6)^T\allowbreak/
(\sqrt{38}/6)
\in E_1 \cap \bar{C}_{1\cdots 5}$, 
which means that the cardinality of the 
orbit $\dS_5 E_1
$ 
cannot exceed $5!/(3!)=20: 
|\dS_5 E_1|\le 20$. 
As for $-\rho R_5$ in \eqref{eq:V2capC1-5},  
the chambers 
$D\subset \calV_2 \ 
(D \in \Dch(\calA^0_5))$ such that 
$\emptyset\ne 
D\cap \bar{C}_{1\cdots 5}\subseteq -\rho R_5$ 
are exactly those five 
chambers given as $D=-\rho E_i, \ 
i=1,\ldots,5$. 


\begin{table}
\caption{$D\cap \bar{C}_{1\cdots 5}\subseteq R_5$ 
($v_1 \ge v_2 \ge v_3 > 0 > v_4 \ge v_5$). 
$\epsilon>0$ is 
small enough.}
\begin{center} 
\begin{tabular}{llr} \toprule 
Defining Inequalities & 
Representative Point & 
$|\dS_5D|$ 
\\ \midrule 
$v_{\{ 3,4\}}>0$ & 
$(
\frac{1}{3}, \ \frac{1}{3}, \ \frac{1}{3}, \ 
-\frac{1}{6}, \ -\frac{5}{6}
)$ &  
20 
\\ 
$v_{\{ 3,4\}}<0, \ v_{\{ 2,4\}}>0$ & 
$(
\frac{1}{3}+\epsilon, \ \frac{1}{3}+\epsilon, \ 
\frac{1}{3}-2\epsilon, \ 
-\frac{1}{3}, \ -\frac{2}{3}
)$ &  
60 
\\ 
$v_{\{ 2,4\}}<0, 
\ v_{\{ 1,5\}}<0, 
\ v_{\{ 1,4\}}>0$ & 
$(
\frac{1}{3}+2\epsilon, \ \frac{1}{3}-\epsilon, \ 
\frac{1}{3}-\epsilon, \ 
-\frac{1}{3}, \ -\frac{2}{3}
)$ &  
60 
\\ 
$v_{\{ 1,5\}}>0 
$ & 
$(
\frac{2}{3}, \ \frac{1}{6}, \ \frac{1}{6}, \ 
-\frac{1}{2}, \ -\frac{1}{2}
)$ &  
30 
\\ 
$v_{\{ 1,4\}}<0
$ &  
$(
\frac{1}{3}, \ \frac{1}{3}, \ \frac{1}{3}, \ 
-\frac{1}{2}, \ -\frac{1}{2}
)$ &  
10 
\\ \bottomrule 
\end{tabular}
\end{center} 
\label{table:R5}
\end{table}


The discussions above imply that 
\begin{eqnarray*}
\calV_2 
\cap \bar{C}_{1\cdots 5} 
&=&
(E_1\cap \bar{C}_{1\cdots 5})
\sqcup \cdots \sqcup 
(E_5\cap \bar{C}_{1\cdots 5}) \\ 
&& \qquad \qquad \sqcup 
((-\rho E_1)\cap \bar{C}_{1\cdots 5})
\sqcup \cdots \sqcup  
((-\rho E_5)\cap \bar{C}_{1\cdots 5}),  
\end{eqnarray*}
and this yields 
$\Dch_2^{1\cdots 5}(\calA_5^0)
=\{ E_1,\ldots,E_5, 
-\rho E_1,\ldots,-\rho E_5\}$ and  
\begin{equation}
\label{eq:c.s.m=5}
\Dch^{1\cdots 5}(\calA_5^0)
=\{ D_1, -\rho D_1, 
E_1,\ldots,E_5, 
-\rho E_1,\ldots,-\rho E_5\}. 
\end{equation}
Since 
$\Dch^{1\cdots 5}(\calA_5^0)$ 
includes 
a cross section 
under the action of $\dS_5$ on $\Dch(\calA_5^0)$, 
we know 
\begin{eqnarray}
\Dch(\calA^0_5)
&=& 
\dS_5D_1 \cup \dS_5(-\rho D_1) 
\cup 
\dS_5E_1 \cup \cdots \cup \dS_5 E_5 \notag \\ 
&& \qquad \qquad \qquad \qquad \quad  
\cup 
\dS_5(-\rho E_1) \cup \cdots \cup 
\dS_5(-\rho E_5), \label{eq:D5=cupcupcup} 
\end{eqnarray} 
which implies 
\begin{eqnarray}
|\Dch(\calA^0_5)|
&\le&  
|\dS_5D_1| + |\dS_5(-\rho D_1)| 
+ \sum_{i=1}^{5}|\dS_5E_i| 
+ \sum_{i=1}^{5}|\dS_5(-\rho E_i)| 
\notag \\ 
&=& 2(|\dS_5D_1|+ \sum_{i=1}^{5}|\dS_5E_i|) 
=2(5+\sum_{i=1}^{5}|\dS_5E_i|).  \label{eq:||<||+||} 
\end{eqnarray} 

From 
the last 
column of Table \ref{table:R5},  
wee can see $\sum_{i=1}^{5}|\dS_5E_i|\le 
20+60+60+30+10=
180$, 
so $|\Dch(\calA^0_5)|\le 
2(5+180)=370$ by \eqref{eq:||<||+||}. 
But since $370$ is equal to
$|\Dch(\calA^0_5)|=|\ch(\calA^0_5)|$ 
(see Lemma \ref{lm:chi(3-8)}), 
the inequality in 
\eqref{eq:||<||+||} is actually an 
equality. 
This means the 
$|\Dch^{1\cdots 5}(\calA_5^0)|=2(1+5)=12$ 
orbits on the right-hand side of 
\eqref{eq:D5=cupcupcup} are all distinct: 
\begin{eqnarray*}
\Dch(\calA^0_5)
&=& 
\dS_5D_1 \sqcup \dS_5(-\rho D_1) 
\sqcup 
\dS_5E_1 \sqcup \cdots \sqcup \dS_5 E_5 \\ 
&& \qquad \qquad \qquad \qquad \quad  
\sqcup 
\dS_5(-\rho E_1) \sqcup \cdots \sqcup 
\dS_5(-\rho E_5).  
\end{eqnarray*} 
Hence, 
$\Dch^{1\cdots 5}(\calA_5^0)$ in 
\eqref{eq:c.s.m=5} is a cross section.  
Therefore, $|\Dch (\calA^0_5)/\dS_5|
=|\Dch^{1\cdots 5}(\calA_5^0)|
=12$, 
and the number of inequivalent ranking patterns of 
unfolding models of codimension one 
is $12-1=11=1+2\cdot 5$. 

Notice, in passing,   
$\sum_{i=1}^{5}|\dS_5E_i|=180$ 
so that the upper bounds in the last 
column of Table \ref{table:R5} 
are actually exact numbers. 


\subsubsection{Case $m=6$} 


When $m=6$, we have 
\[
\calV_2\cap \bar{C}_{1\cdots 6}
=\bigsqcup_{
D\subset \calV_2, \ 
D \in \Dch(\calA^0_6)
}
(D\cap \bar{C}_{1\cdots 6})
=R_{6,1}\sqcup(-\rho R_{6,1})
\sqcup R_{6,2}\sqcup (-\rho R_{6,2}), 
\] 
where  
\begin{eqnarray*}
R_{6,1}
&:=& \{ (v_1,\ldots,v_6)^T\in \bbS^4: 
v_1\ge v_2 \ge v_3 \ge v_4 > 0 > v_5 \ge v_6\}, \\ 
R_{6,2}
&:=& \{ (v_1,\ldots,v_6)^T\in \bbS^4: 
v_1\ge v_2 \ge v_3 > 0 > v_4 \ge v_5 \ge v_6, \ v_3>-v_4\}.  
\end{eqnarray*}



\begin{table}
\caption{$D\cap \bar{C}_{1\cdots 6} \subseteq R_{6,1}$ 
($v_1 \ge v_2 \ge v_3 \ge v_4 >0>v_5 \ge v_6$). 
$\epsilon>0$ is small enough.}
\begin{center} 
\setlength{\tabcolsep}{2pt}
\begin{tabular}{llr} \toprule 
{\footnotesize 
Defining Inequalities} & 
{\footnotesize 
Representative Point} & 
{\footnotesize 
$|\dS_6 D|$}  
\\ \midrule 
{\footnotesize 
$v_{\{ 4, 5\}}>0$} & 
{\footnotesize 
$(
\frac{1}{4}, \ 
\frac{1}{4}, \ 
\frac{1}{4}, \ 
\frac{1}{4}, \ 
-\frac{1}{4}+\epsilon, \ 
-\frac{3}{4}-\epsilon
)$} & 
{\footnotesize 
30} 
\\ 
{\footnotesize 
$v_{\{ 4, 5\}}<0, \ v_{\{ 3, 5\}}>0$} & 
{\footnotesize 
$(
\frac{1}{4}+\epsilon, \ 
\frac{1}{4}+\epsilon, \ 
\frac{1}{4}+\epsilon, \ 
\frac{1}{4}-3\epsilon, \ 
-\frac{1}{4}, \ 
-\frac{3}{4}
)$} & 
{\footnotesize 
120} 
\\ 
{\footnotesize 
$v_{\{ 3, 5\}}<0, \ v_{\{ 2, 5\}}>0, \ v_{\{ 3, 4, 5\}}>0$} &  
{\footnotesize 
$(
\frac{1}{4}+\epsilon, \ 
\frac{1}{4}+\epsilon, \ 
\frac{1}{4}-\epsilon, \ 
\frac{1}{4}-\epsilon, \ 
-\frac{1}{4}, \ 
-\frac{3}{4}
)$} & 
{\footnotesize 
180} 
\\ 
{\footnotesize 
$v_{\{ 2, 5\}}<0, \ v_{\{ 3, 4, 5\}}>0, \ v_{\{ 1, 5\}}>0$} & 
{\footnotesize  
$(
\frac{3}{4}, \ 
\frac{1}{12}, \ 
\frac{1}{12}, \ 
\frac{1}{12}, \ 
-\frac{1}{8}, \ 
-\frac{7}{8}
)$} & 
{\footnotesize  
120} 
\\ 
{\footnotesize 
$v_{\{ 2, 5\}}<0, 
\ v_{\{ 3, 4, 5\}}<0, \ v_{\{ 2, 4, 5\}}>0, 
\ v_{\{ 1, 5\}}>0$} &  
{\footnotesize 
$(
\frac{3}{4}-\epsilon, \ 
\frac{1}{12}+\epsilon, \ 
\frac{1}{12}, \ 
\frac{1}{12}, \ 
-\frac{1}{6}-\frac{\epsilon}{2}, \ 
-\frac{5}{6}+\frac{\epsilon}{2}
)$} &  
{\footnotesize 
360} 
\\ 
{\footnotesize 
$v_{\{ 2, 4, 5\}}<0, \ v_{\{ 2, 3, 5\}}>0, \ v_{\{ 1, 5\}}>0$} & 
{\footnotesize  
$(
\frac{3}{4}+\epsilon, \ 
\frac{1}{12}, \ 
\frac{1}{12}, \ 
\frac{1}{12}-\epsilon, \ 
-\frac{1}{6}+\frac{\epsilon}{2}, \ 
-\frac{5}{6}-\frac{\epsilon}{2}
)$} &  
{\footnotesize 
360} 
\\ 
{\footnotesize 
$v_{\{ 2, 3, 5\}}<0, 
\ v_{\{ 1, 6\}}<0, 
\ v_{\{ 1, 5\}}>0$} &  
{\footnotesize 
$(
\frac{3}{4}, \ 
\frac{1}{12}, \ 
\frac{1}{12}, \ 
\frac{1}{12}, \ 
-\frac{1}{5}, \ 
-\frac{4}{5}
)$} &  
{\footnotesize 
120} 
\\ 
{\footnotesize 
$
v_{\{ 1, 6 \}}>0 
$} & 
{\footnotesize 
$(
\frac{3}{4}, \ 
\frac{1}{12}, \ 
\frac{1}{12}, \ 
\frac{1}{12}, \ 
-\frac{1}{2}, \ 
-\frac{1}{2}
)$} &  
{\footnotesize 
60} 
\\ 
{\footnotesize 
$v_{\{ 3, 4, 5\}}<0, \ v_{\{ 2, 5\}}>0$} & 
{\footnotesize 
$(
\frac{5}{12}, \ 
\frac{5}{12}, \ 
\frac{1}{12}, \ 
\frac{1}{12}, \ 
-\frac{1}{4}, \ 
-\frac{3}{4} 
)$} &  
{\footnotesize 
180} 
\\ 
{\footnotesize 
$v_{\{ 1, 5\}}<0, \ v_{\{ 3, 4, 5\}}>0$} & 
{\footnotesize  
$(
\frac{1}{4}, \ 
\frac{1}{4}, \ 
\frac{1}{4}, \ 
\frac{1}{4}, \ 
-\frac{1}{2}+\epsilon, \ 
-\frac{1}{2}-\epsilon 
)$} & 
{\footnotesize  
30} 
\\ 
{\footnotesize 
$v_{\{ 1, 5\}}<0, 
\ v_{\{ 3, 4, 5\}}<0, \ v_{\{ 2, 4, 5\}}>0$} &  
{\footnotesize 
$(
\frac{1}{4}+\epsilon, \ 
\frac{1}{4}+\epsilon, \ 
\frac{1}{4}-\epsilon, \ 
\frac{1}{4}-\epsilon, \ 
-\frac{1}{2}+\epsilon, \ 
-\frac{1}{2}-\epsilon
)$} & 
{\footnotesize 
180} 
\\ 
{\footnotesize 
$v_{\{ 1, 5\}}<0, 
\ v_{\{ 2, 4, 5\}}<0, 
\ v_{\{ 2, 3, 5\}}>0, 
\ v_{\{ 1, 4, 5\}}>0$} &  
{\footnotesize 
$(
\frac{1}{4}+\epsilon, \ 
\frac{1}{4}, \ 
\frac{1}{4}, \ 
\frac{1}{4}-\epsilon, \ 
-\frac{1}{2}+\frac{\epsilon}{2}, \ 
-\frac{1}{2}-\frac{\epsilon}{2} 
)$} &  
{\footnotesize 
360} 
\\ 
{\footnotesize 
$v_{\{ 1, 4, 5\}}<0$} &  
{\footnotesize 
$(
\frac{1}{4}+\epsilon, \ 
\frac{1}{4}+\epsilon, \ 
\frac{1}{4}+\epsilon, \ 
\frac{1}{4}-3\epsilon, \ 
-\frac{1}{2}, \ 
-\frac{1}{2} 
)$} &  
{\footnotesize 
60} 
\\ 
{\footnotesize 
$v_{\{ 2, 3, 5\}}<0, \ v_{\{ 1, 5\}}<0$} & 
{\footnotesize  
$(
\frac{1}{4}+3\epsilon, \ 
\frac{1}{4}-\epsilon, \ 
\frac{1}{4}-\epsilon, \ 
\frac{1}{4}-\epsilon, \ 
-\frac{1}{2}, \ 
-\frac{1}{2}
)$} & 
{\footnotesize  
60} 
\\ \bottomrule 
\end{tabular}
\end{center}
\label{table:R6-1}
\end{table}


\begin{table}
%
\caption{$D\cap \bar{C}_{1\cdots 6} \subseteq R_{6,2}$ 
($v_1 \ge v_2 \ge v_3 > 0 > v_4 \ge v_5 \ge v_6, \ 
v_3 > -v_4$). 
$\epsilon>0$ is small enough.}
\begin{center} 
\setlength{\tabcolsep}{2pt}
\begin{tabular}{llr} \toprule 
{\footnotesize 
Defining Inequalities} & 
{\footnotesize 
Representative Point} & 
{\footnotesize 
$|\dS_6 D|$}  
\\ \midrule 
{\footnotesize 
$v_{\{ 3,4,5\}}>0$} &  
{\footnotesize 
$(
\frac{1}{3}, \ \frac{1}{3}, \ \frac{1}{3}, \ 
-\frac{1}{6}+\epsilon, \ -\frac{1}{6}+\epsilon, \ 
-\frac{2}{3}-2\epsilon 
)$} &  
{\footnotesize 
60} 
\\ 
{\footnotesize 
$v_{\{ 3,5\}}>0, \ v_{\{ 3,4,5\}}<0, \ v_{\{ 2,4,5\}}>0$} & 
{\footnotesize 
$(
\frac{1}{3}+\epsilon, \ \frac{1}{3}+\epsilon, \ 
\frac{1}{3}-2\epsilon, \ 
-\frac{1}{6}, \ -\frac{1}{6}, \ 
-\frac{2}{3} 
)$} &  
{\footnotesize 
180} 
\\ 
{\footnotesize 
$v_{\{ 3,5\}}>0, \ v_{\{ 2,4,5\}}<0, 
\ v_{\{ 1,6 \}}<0, 
\ v_{\{ 1,4,5\}}>0$} &  
{\footnotesize 
$(
\frac{2}{3}, \ 
\frac{1}{6}, \ \frac{1}{6}, \ 
-\frac{1}{8}, \ -\frac{1}{8}, \ 
-\frac{3}{4} 
)$} &  
{\footnotesize 
180} 
\\ 
{\footnotesize 
$
v_{\{ 3,5\}}<0, \ v_{\{ 2,4,5\}}>0$} &  
{\footnotesize 
$(
\frac{1}{3}+\epsilon, \ \frac{1}{3}+\epsilon, \ 
\frac{1}{3}-2\epsilon, \ 
-\epsilon, \ 
-\frac{1}{3}+\epsilon, \ 
-\frac{2}{3} 
)$} &  
{\footnotesize 
360} 
\\ 
{\scriptsize
$
v_{\{ 3,5\}}<0,  v_{\{ 2,5\}}>0, 
 v_{\{ 2,4,5\}}<0, 
 v_{\{ 1,6\}}<0, 
 v_{\{ 1,4,5\}}>0$} & 
{\footnotesize  
$(
\frac{2}{3}, \ 
\frac{1}{6}+\epsilon, \ 
\frac{1}{6}-\epsilon, \ 
-2\epsilon, \ 
-\frac{1}{6}, \ 
-\frac{5}{6}+2\epsilon 
)$} &  
{\footnotesize 
720} 
\\ 
{\footnotesize 
$
v_{\{ 2,5\}}<0, 
\ v_{\{ 1,6\}}<0, 
\ v_{\{ 1,4,5\}}>0$} & 
{\footnotesize  
$(
\frac{2}{3}, \ 
\frac{1}{6}, \ \frac{1}{6}, \ 
-\epsilon, \ 
-\frac{1}{4}, \ 
-\frac{3}{4}+\epsilon 
)$} & 
{\footnotesize  
360} 
\\ 
{\footnotesize 
$
v_{\{ 2,3,5\}}>0, 
\ v_{\{ 1,6\}}>0, 
\ v_{\{ 2,3,6\}}<0$} & 
{\footnotesize 
$(
\frac{2}{3}, \ 
\frac{1}{6}, \ \frac{1}{6}, \ 
-2\epsilon, \ 
-\frac{1}{3}+\epsilon, \ 
-\frac{2}{3}+\epsilon 
)$} &  
{\footnotesize 
360} 
\\ 
{\footnotesize 
$
v_{\{ 2,3,6\}}>0, 
\ v_{\{ 1,6\}}>0$} &  
{\footnotesize 
$(
\frac{1}{2}, \ 
\frac{1}{4}, \ \frac{1}{4}, \ 
-\epsilon, \ 
-\frac{1}{2}+\frac{\epsilon}{2}, \ 
-\frac{1}{2}+\frac{\epsilon}{2} 
)$} &  
{\footnotesize 
180} 
\\ 
{\footnotesize 
$
v_{\{ 2,3,5\}}<0$} &  
{\footnotesize 
$(
\frac{2}{3}, \  
\frac{1}{6}, \ \frac{1}{6}, \ 
-\epsilon, \ 
-\frac{1}{2}+\frac{\epsilon}{2}, \ 
-\frac{1}{2}+\frac{\epsilon}{2} 
)$} & 
{\footnotesize  
180} 
\\ 
{\footnotesize 
$v_{\{ 3,5\}}>0, \ v_{\{ 1,4,5\}}<0, \ v_{\{ 1,6\}}<0$} & 
{\footnotesize  
$(
\frac{1}{3}, \ \frac{1}{3}, \ \frac{1}{3}, \ 
-\frac{1}{6}-\epsilon, \ 
-\frac{1}{6}-\epsilon, \ 
-\frac{2}{3}+2\epsilon 
)$} &  
{\footnotesize 
60} 
\\ 
{\footnotesize 
$
v_{\{ 3,5\}}<0, \ v_{\{ 2,5\}}>0, 
\ v_{\{ 1,4,5\}}<0, \ v_{\{ 1,6\}}<0$} & 
{\footnotesize 
$(
\frac{1}{3}+\epsilon, \ 
\frac{1}{3}+\epsilon, \ 
\frac{1}{3}-2\epsilon, \ 
-3\epsilon, \ 
-\frac{1}{3}+\epsilon, \ 
-\frac{2}{3}+2\epsilon 
)$} & 
{\footnotesize  
360} 
\\ 
{\footnotesize 
$
v_{\{ 2,5\}}<0, \ v_{\{ 1,5\}}>0, 
\ v_{\{ 1,4,5\}}<0, \ v_{\{ 1,6\}}<0$} & 
{\footnotesize 
$(
\frac{1}{3}+2\epsilon, \ 
\frac{1}{3}-\epsilon, \ 
\frac{1}{3}-\epsilon, \ 
-3\epsilon, \ 
-\frac{1}{3}, \ 
-\frac{2}{3}+3\epsilon 
)$} &  
{\footnotesize 
360} 
\\ 
{\footnotesize 
$
v_{\{ 1,5\}}<0$} &  
{\footnotesize 
$(
\frac{1}{3}, \ \frac{1}{3}, \ \frac{1}{3}, \ 
-2\epsilon, \ 
-\frac{1}{2}+\epsilon, \ 
-\frac{1}{2}+\epsilon 
)$} &  
{\footnotesize 
60} 
\\ \bottomrule 
\end{tabular} 
\end{center} 
\label{table:R6-2}
\end{table}


All nonempty 
$D\cap \bar{C}_{1\cdots 6}$ 
($D\subset \calV_2, \ 
D \in \Dch(\calA^0_6)$) 
that are included 
in $R_{6,1}$ 
and in $R_{6,2}$ 
are listed in Tables \ref{table:R6-1} and 
\ref{table:R6-2}, respectively. 
The first columns 
provide 
the defining inequalities of each 
$D\cap \bar{C}_{1\cdots 6}$ 
(besides $v_1\ge v_2 \ge v_3 \ge v_4 > 0 > v_5 \ge v_6$ 
and $v_1\ge v_2 \ge v_3 > 0 > v_4 \ge v_5 \ge v_6, 
\ v_3>-v_4$, respectively); 
the second columns 
show 
an unnormalized 
representative point of  
each $D\cap \bar{C}_{1\cdots 6}$; and 
the last columns 
include 
an upper bound for 
$|\dS_6D|$ for each $D$. 



There are $14$ (resp. $13$) 
rows 
in Table \ref{table:R6-1} 
(resp. Table \ref{table:R6-2}), 
and the sum of the upper bounds for 
$|\dS_6D|$ 
in the last column of the table 
is $2220$ (resp. $3420$). 
Since the value $2(6+2220+3420)=11292$ 
equals 
$|\ch(\calA^0_6)|$ (Lemma \ref{lm:chi(3-8)}), 
we can conclude, by the same reasoning as in 
the case of $m=5$, that 
the number of 
inequivalent ranking patterns of unfolding models of codimension one 
is $1+2(14+13)=55$. 


\

\noindent 
{\bf Open problem}: 
We have seen that for any $m\le 6$, 
subset $\Dch^{1\cdots m}(\calA^0)\subset \Dch(\calA^0)$ 
is a cross section so that the upper bound in  
\eqref{eq:1+||} is actually the exact number. 
Does this hold true for all $m$? 

\ 

\ 

\noindent
{\bf
Acknowledgments.} 
The authors are very grateful to two anonymous referees 
for their valuable suggestions that helped to improve 
the presentation of an earlier version of this paper.

\end{document}